\documentclass{amsart}
\usepackage{amssymb}

\numberwithin{equation}{section}
\input{tcilatex}

\begin{document}
\title[Uniqueness of minimizers]{Uniqueness of generalized $p$-area
minimizers and integrability of a horizontal normal in the Heisenberg group}
\author{Jih-Hsin Cheng}
\address{Institute of Mathematics, Academia Sinica, Taipei and National
Center for Theoretical Sciences, Taipei Office, Taiwan, R.O.C.}
\email{ cheng@math.sinica.edu.tw}
\urladdr{http://www.math.sinica.edu.tw}
\thanks{}
\author{Jenn-Fang Hwang}
\address{Institute of Mathematics, Academia Sinica, Taipei, Taiwan, R.O.C.}
\email{ majfh@math.sinica.edu.tw}
\urladdr{http://www.math.sinica.edu.tw}
\thanks{}
\thanks{}

\begin{abstract}
We study the uniqueness of generalized $p$-minimal surfaces in the
Heisenberg group. The generalized $p$-area of a graph defined by $u$ reads $%
\int |\nabla u+\vec{F}|+Hu.$ If $u$ and $v$ are two minimizers for the
generalized $p$-area satisfying the same Dirichlet boundary condition, then
we can only get $N_{\vec{F}}(u)$ $=$ $N_{\vec{F}}(v)$ (on the nonsingular
set) where $N_{\vec{F}}(w)$ $:=$ $\frac{\nabla w+\vec{F}}{|\nabla w+\vec{F}|}%
.$ To conclude $u$ $=$ $v$ (or $\nabla u$ $=$ $\nabla v)$, it is not
straightforward as in the Riemannian case, but requires some special
argument in general. In this paper, we prove that $N_{\vec{F}}(u)$ $=$ $N_{%
\vec{F}}(v)$ implies $\nabla u$ $=$ $\nabla v$ in dimension $\geq $ $3$
under some rank condition on derivatives of $\vec{F}$ or the
nonintegrability condition of contact form associated to $u$ or $v$. Note
that in dimension 2 ($n=1),$ the above statement is no longer true. Inspired
by an equation for the horizontal normal $N_{\vec{F}}(u),$ we study the
integrability for a unit vector to be the horizontal normal of a graph. We
find a Codazzi-like equation together with this equation to form an
integrability condition.
\end{abstract}

\maketitle


\bigskip



\section{\textbf{Introduction and statement of the results}}

Recall that the $p$-area (pseudohermitian area or called horizontal area by
some authors) is a special case of the generalized $p$-area:%
\begin{equation}
\mathcal{F}_{H}(u)\equiv \int_{\Omega }\{|\nabla u+\vec{F}|+Hu\}\text{ }%
dx^{1}\wedge dx^{2}\wedge ...\wedge dx^{m}.  \label{eqn1.0}
\end{equation}

\noindent where $\Omega \subset R^{m}$ is a bounded domain, $u$ $\in $ $%
W^{1,1}(\Omega ),$ $\vec{F}$ is an $L^{1}$ vector field on $\Omega ,$ and $H$
$\in $ $L^{\infty }(\Omega ),$ say. We denote $\mathcal{F}_{H}$ by $\mathcal{%
F}_{0}$ for the case of $H$ $=$ $0:$%
\begin{equation}
\mathcal{F}_{0}(u)\equiv \int_{\Omega }|\nabla u+\vec{F}|.  \label{eqn1.1}
\end{equation}%
\noindent $\mathcal{F}_{0}\mathcal{(\cdot )}$ is called the $p$-area (of the
graph defined by $u$ over $\Omega )$ if $\vec{F}$ $=$ $-\vec{X}^{\ast }$
where $\vec{X}^{\ast }$ $=$ $(x^{1^{\prime }},$ $-x^{1},$ $x^{2^{\prime }},$ 
$-x^{2},...,$ $x^{n^{\prime }},$ $-x^{n}),$ $m$ $=$ $2n$ (see \cite{chmy}).
In the case of a graph $\Sigma $ over the $R^{2n}$-hyperplane in the
Heisenberg group, the above definition of $p$-area coincides with those
given in \cite{CDG}, \cite{DGN}, and \cite{Pau}. In particular these
notions, especially in the framework of geometric measure theory, have been
used to study existence or regularity properties of minimizers for the
relative perimeter or extremizers of isoperimetric inequalities (see, e.g., 
\cite{DGN}, \cite{GN}, \cite{LM}, \cite{LR}, \cite{MR}, \cite{Pan}, \cite{RR}%
).

The $p$-area can also be identified with the $2n+1$-dimensional spherical
Hausdorff measure of $\Sigma $ (see, e.g., \cite{Ba}, \cite{FSS}). Some
authors take the viewpoint of so called intrinsic graphs (see, e.g., \cite%
{FSS}, \cite{ASCV}, \cite{BSC}). Starting from the work \cite{chmy} (see
also \cite{ch}), we studied the subject from the viewpoint of partial
differential equations and that of differential geometry (see \cite{chy}, 
\cite{chy1}, \cite{chmy2}, and \cite{ch1}; we use the term $p$-minimal since
this is the notion of minimal surfaces in pseudohermitian geometry; "$p$"
stands for "pseudohermitian"). In particular, the generalized $p$-area (\ref%
{eqn1.0}) has been studied in the Heisenberg group.

First look at the integrand $D_{u}$ $:=$ $|\nabla u+\vec{F}|$ in $\mathcal{F}%
_{H}.$ Denote $\frac{\partial u}{\partial x_{i}}$ by $u_{i}.$ We compute%
\begin{equation*}
\frac{\partial ^{2}D_{u}}{\partial u_{i}\partial u_{j}}=\frac{\delta _{ij}}{%
D_{u}}-\frac{(u_{i}+F_{i})(u_{j}+F_{j})}{D_{u}^{3}}.
\end{equation*}%
\noindent Observe that $\frac{\partial ^{2}D_{u}}{\partial u_{i}\partial
u_{j}}\xi _{i}\xi _{j}$ $\geq $ $0$ (summation convention), but $\frac{%
\partial ^{2}D_{u}}{\partial u_{i}\partial u_{j}}\xi _{i}\xi _{j}$ $=$ $0$
does not imply $\xi _{i}$ $=$ $0$ for all $i.$ ($\frac{\partial ^{2}D_{u}}{%
\partial u_{i}\partial u_{j}})$ being not positive definite causes trouble
in studying $\mathcal{F}_{H}.$ Let $S_{\vec{F}}(w)$ denote the singular set
of a real function $w$ defined on $\Omega ,$ which consists of points $p$ $%
\in $ $\Omega $ such that $\nabla w+\vec{F}$ $=$ $0$ at $p.$ Let $\mathcal{%
F(\varepsilon )}$ $:=$ $\mathcal{F}_{H}(u_{\varepsilon })\mathcal{\ }$where $%
u_{\varepsilon }$ $:=$ $u$ $+$ $\varepsilon \varphi $, $\varphi $ $\in $ $%
W_{0}^{1,1}.$ We have%
\begin{eqnarray}
\frac{d\mathcal{F(}0\pm \mathcal{)}}{d\varepsilon } &:&=\lim_{\varepsilon
\rightarrow 0\pm }\frac{\mathcal{F(\varepsilon )-F(}0)}{\varepsilon }
\label{1.1.3} \\
&=&\pm \int_{S_{\vec{F}}(u)}|\nabla \varphi |\text{ }+\text{ }\int_{\Omega
\backslash S_{\vec{F}}(u)}N_{\vec{F}}(u)\cdot \nabla \varphi +\int_{\Omega
}H\varphi  \notag
\end{eqnarray}

\noindent ((3.3) with $\hat{\varepsilon}$ $=$ $0$ in \cite{chy}) where we
denote the horizontal normal $\frac{\nabla w+\vec{F}}{|\nabla w+\vec{F}|}$
for a real function $w$ ($\in $ $W^{1}(\Omega ),$ say) by $N_{\vec{F}}(w)$
(or $\nu ^{w}$; the notation $N_{\vec{F}}(w)$ has been used previously. But
the notation $\nu ^{w}$ is concise)$.$ From (\ref{1.1.3}), the first
variation formula of $\mathcal{F}_{H}$ at $u,$ we found that $\pm \int_{S_{%
\vec{F}}(u)}|\nabla \varphi |$ is not negligible if $H_{m}(S_{\vec{F}}(u))$ $%
\neq $ $0$. In \cite{ch1}, we extended the range of $u$ $\in $ $W^{1,1}$ to $%
u$ $\in $ $BV$ (see also \cite{SCV}, \cite{PCTV}) and computed the first and
second variations. For $u,v$ $\in $ $BV,$ $u_{\varepsilon }$ $=$ $u$ $+$ $%
\varepsilon \varphi $ with $\varphi $ $=$ $v$ $-$ $u,$ the right and left
derivatives $\mathcal{F}_{\pm }^{\prime }\mathcal{(\varepsilon )}$ $:=\lim_{%
\tilde{\varepsilon}\rightarrow \varepsilon \pm }\frac{\mathcal{F(\tilde{%
\varepsilon})-F(}\varepsilon )}{\tilde{\varepsilon}-\varepsilon }$ exist and
satisfy%
\begin{equation*}
\mathcal{F}_{-}^{\prime }\mathcal{(\varepsilon }_{1}\mathcal{)\leq F}%
_{+}^{\prime }\mathcal{(\varepsilon }_{1}\mathcal{)\leq F}_{-}^{\prime }%
\mathcal{(\varepsilon }_{2}\mathcal{)\leq F}_{+}^{\prime }\mathcal{%
(\varepsilon }_{2}\mathcal{)}
\end{equation*}

\noindent for $\varepsilon _{1}$ $<$ $\varepsilon _{2}$ ((3.20) in \cite{ch1}%
). We also have%
\begin{equation*}
\lim_{\mathcal{\varepsilon }_{2}\rightarrow \mathcal{\varepsilon }_{1}+}%
\mathcal{F}_{\pm }^{\prime }\mathcal{(\varepsilon }_{2}\mathcal{)=F}%
_{+}^{\prime }\mathcal{(\varepsilon }_{1}\mathcal{)}\text{, }\lim_{\mathcal{%
\varepsilon }_{2}\rightarrow \mathcal{\varepsilon }_{1}-}\mathcal{F}_{\pm
}^{\prime }\mathcal{(\varepsilon }_{2}\mathcal{)=F}_{-}^{\prime }\mathcal{%
(\varepsilon }_{1}\mathcal{)}
\end{equation*}

\noindent and $\mathcal{F}$ is convex in $\varepsilon $ ((3.21) in \cite{ch1}%
). For the second variation, although $\mathcal{F}_{+}^{\prime }\mathcal{%
(\varepsilon )}$ may not equal $\mathcal{F}_{-}^{\prime }\mathcal{%
(\varepsilon )}$, we have%
\begin{equation*}
\lim_{\mathcal{\varepsilon }_{2}\rightarrow \mathcal{\varepsilon }_{1}+}%
\frac{\mathcal{F}_{\pm }^{\prime }\mathcal{(\varepsilon }_{2}\mathcal{)-F}%
_{+}^{\prime }\mathcal{(\varepsilon }_{1}\mathcal{)}}{\mathcal{\varepsilon }%
_{2}-\mathcal{\varepsilon }_{1}}=\lim_{\mathcal{\varepsilon }_{2}\rightarrow 
\mathcal{\varepsilon }_{1}-}\frac{\mathcal{F}_{\pm }^{\prime }\mathcal{%
(\varepsilon }_{2}\mathcal{)-F}_{-}^{\prime }\mathcal{(\varepsilon }_{1}%
\mathcal{)}}{\mathcal{\varepsilon }_{2}-\mathcal{\varepsilon }_{1}}
\end{equation*}

\noindent (Theorem C in \cite{ch1}). That is to say, the first variation may
have jumps, but the right and left limits of the second variation exist and
coincide. This is an interesting property.

In \cite{chy}, we proved the uniqueness of minimizers for the generalized $p$%
-area $\mathcal{F}_{H}$ in the space $W^{1,2}$ among other things. Recall
that $u$ $\in $ $W^{1,1}(\Omega )$ ($W^{1,2}(\Omega ),$ resp.) is called a
minimizer for $\mathcal{F}_{H}$ (see (\ref{eqn1.0})) if there holds 
\begin{equation*}
\mathcal{F}_{H}(u)\leq \mathcal{F}_{H}(u+\varphi )
\end{equation*}%
\noindent for any $\varphi $ $\in $ $W_{0}^{1,1}(\Omega )$ ($%
W_{0}^{1,2}(\Omega ),$ resp.). Let $\vec{F}^{\ast }$ $:=$ $(F_{2},$ $-F_{1},$
$F_{4},$ $-F_{3},...,$ $F_{2n},$ $-F_{2n-1})$ for $\vec{F}$ $=$ $(F_{1},$ $%
F_{2},...,$ $F_{2n}).$

\bigskip

\textbf{Theorem 1.1 (Theorem B in \cite{chy})}. \textit{Let }$\Omega $%
\textit{\ be a bounded domain in }$R^{2n}.$\textit{\ Let }$u,v\in
W^{1,2}(\Omega )$\textit{\ be two minimizers for }$\mathcal{F}_{H}$\textit{\
such that }%
\begin{equation}
u-v\in W_{0}^{1,2}(\Omega ).  \label{1.1.1}
\end{equation}

\textit{\noindent Suppose }$H$\textit{\ }$\in $\textit{\ }$L^{\infty
}(\Omega )$\textit{\ and }$\vec{F}$\textit{\ }$\in $\textit{\ }$%
W^{1,2}(\Omega )$\textit{\ satisfying}%
\begin{equation}
\func{div}\mathit{\vec{F}}^{\ast }\mathit{>0}\text{ a.e. (}\func{div}\mathit{%
\vec{F}}^{\ast }\mathit{<0}\text{ a.e., resp.)}  \label{1.1.2}
\end{equation}

\textit{\noindent Then }$u\equiv v$\textit{\ in }$\Omega $\textit{\ (a.e.).}

\bigskip

The uniqueness of $BV$ solutions to the appropriate Dirichlet problem is
still unknown. However for $u,v$ $\in $ $W^{1,2}$ as in Theorem 1.1 (Theorem
B in \cite{chy}), since $\mathcal{F(\varepsilon )}$ is nondecreasing and $%
\mathcal{F(}0\mathcal{)}$ $=$ $\mathcal{F}_{H}(u)$ $=$ $\mathcal{F}_{H}(v)$ $%
=$ $\mathcal{F(}1\mathcal{)}$ ($u,v$ being minimizers for $\mathcal{F}_{H}),$
we have $\mathcal{F(\varepsilon )}$ $=$ $\mathcal{F(}0\mathcal{)}$ $=$ $%
\mathcal{F(}1\mathcal{)}$ for all $\varepsilon ,$ $0$ $\leq $ $\varepsilon $ 
$\leq $ $1.$ Moreover, we can show that there are at most countably many $%
\varepsilon $ such that $H_{m}(S_{\vec{F}}(u_{\varepsilon }))$ $\neq $ $0.$
Choose $\varepsilon _{1},$ $\varepsilon _{2}$ $\in $ $(0,1)$ with $H_{m}(S_{%
\vec{F}}(u_{\varepsilon _{1}}))$ $=$ $H_{m}(S_{\vec{F}}(u_{\varepsilon
_{2}}))$ $=$ $0.$ We then have $\mathcal{F(\varepsilon }_{1}\mathcal{)}$ $=$ 
$\mathcal{F(\varepsilon }_{2}\mathcal{)}$ and $\mathcal{F}^{\prime }\mathcal{%
(\varepsilon }_{1}\mathcal{)}$ $=$ $\mathcal{F}^{\prime }\mathcal{%
(\varepsilon }_{2}\mathcal{)}$ $=$ $0$. It follows from (\ref{1.1.3}) that%
\begin{eqnarray*}
0 &=&\mathcal{F}^{\prime }\mathcal{(\varepsilon }_{2}\mathcal{)}-\mathcal{F}%
^{\prime }\mathcal{(\varepsilon }_{1}\mathcal{)} \\
&=&\int_{\Omega }(N_{\vec{F}}(u_{\varepsilon _{2}})-N_{\vec{F}%
}(u_{\varepsilon _{1}}))\cdot \frac{\nabla u_{\varepsilon _{2}}-\nabla
u_{\varepsilon _{1}}}{\mathcal{\varepsilon }_{2}-\mathcal{\varepsilon }_{1}}.
\end{eqnarray*}

\noindent Here we have used $\varphi $ $=$ $v-u$ $=$ $\frac{u_{\varepsilon
_{2}}-u_{\varepsilon _{1}}}{\mathcal{\varepsilon }_{2}-\mathcal{\varepsilon }%
_{1}}.$ Now the equality $(N_{\vec{F}}(u_{\varepsilon _{2}})-N_{\vec{F}%
}(u_{\varepsilon _{1}}))\cdot (\nabla u_{\varepsilon _{2}}-\nabla
u_{\varepsilon _{1}})$ $=$ $\frac{1}{2}(D_{u_{\varepsilon
_{2}}}+D_{u_{\varepsilon _{1}}})|N_{\vec{F}}(u_{\varepsilon _{2}})-N_{\vec{F}%
}(u_{\varepsilon _{1}})|^{2}$ (cf. Lemma 5.1$^{\prime }$ in \cite{chmy})
forces that $N_{\vec{F}}(u_{\varepsilon _{1}})$ $=$ $N_{\vec{F}%
}(u_{\varepsilon _{2}})$ (a.e.).

Note that having made use of the boundary condition (\ref{1.1.1}), we prove $%
N_{\vec{F}}(u_{\varepsilon _{1}})$ $=$ $N_{\vec{F}}(u_{\varepsilon _{2}})$
for so called regular $\varepsilon _{1},$ $\varepsilon _{2}$ $\in $ $[0,1]$
(see Section 3 or \cite{chy} for more detail)$,$ in which $u_{\varepsilon }$ 
$:=$ $u$ $+$ $\varepsilon (v-u)$ (in the case of good regularity, we have $%
N_{\vec{F}}(u)$ $=$ $N_{\vec{F}}(v))$. In fact, the difficulty of the proof
of Theorem 1.1 is that we may have $H_{2n}(S_{\vec{F}}(u))$ $\neq $ $0$ or $%
H_{2n}(S_{\vec{F}}(v))$ $\neq $ $0.$ We avoid such difficulty by working on
regular $\varepsilon .$ Next together with the condition (\ref{1.1.2}) we
can show $u\equiv v$\textit{\ }(in particular, $\nabla u$ $=$ $\nabla v)$ in 
$\Omega $ (a.e.).

In this paper we will first focus on the problem when $N_{\vec{F}}(u)$ $=$ $%
N_{\vec{F}}(v)$ implies $\nabla u$ $=$ $\nabla v$ with no boundary condition
(\ref{1.1.1}). In general this is not possible. For instance, $u$ $=$ $xy$
and $v$ $=$ $xy$ $+$ $y$ in the Heisenberg group of dimension 3. See Example
2.2 for details. On the positive side, we find a rank condition on the
derivatives of $\vec{F}.$ Let $h_{IJ}$ $:=$ $\partial _{I}F_{J}-\partial
_{J}F_{I}$ (see (\ref{A7})). The rank of a matrix $A,$ denoted as $rank(A),$
is the dimension of the range $Range(A)$ (or image) of $A.$ Note that for
all the results below in this paper, we do not assume $m$ $=$ $2n.$

\bigskip

\textbf{Theorem A.} \textit{Let }$u,$\textit{\ }$v$\textit{\ }$\in $\textit{%
\ }$C^{2}(\Omega )$ \textit{and }$\vec{F}$\textit{\ }$\in $\textit{\ }$%
C^{1}(\Omega )$\textit{\ where }$\Omega $\textit{\ is a domain in }$R^{m}.$%
\textit{\ Suppose both }$N_{\vec{F}}(u)$\textit{\ }$=$\textit{\ }$N_{\vec{F}%
}(v)$\textit{\ and }$\nabla N_{\vec{F}}(u)$ $=$ $\nabla N_{\vec{F}}(v)$ 
\textit{at one point }$p$\textit{\ }$\in $ $\Omega $\textit{\ }$\backslash $%
\textit{\ }$[S_{\vec{F}}(u)$\textit{\ }$\cup $\textit{\ }$S_{\vec{F}}(v)].$%
\textit{\ Assume }%
\begin{equation*}
m\mathit{\ }\geq \mathit{\ }rank(h_{IJ}(p))\mathit{\ }\geq \mathit{\ }3.
\end{equation*}%
\textit{Then }$\nabla u$\textit{\ }$=$\textit{\ }$\nabla v$ \textit{at }$p$%
\textit{\ }$\in $ $\Omega $\textit{\ }$\backslash $\textit{\ }$[S_{\vec{F}%
}(u)$\textit{\ }$\cup $\textit{\ }$S_{\vec{F}}(v)].$

\bigskip

By adding the boundary condition we then have the uniqueness of minimizers
for $\mathcal{F}_{H}.$

\bigskip

\textbf{Corollary A.1}. \textit{Let }$\Omega $\textit{\ be a bounded domain
in }$R^{m}.$\textit{\ Take }$\vec{F}$\textit{\ }$\in $\textit{\ }$C^{1}(\bar{%
\Omega})$\textit{\ and }$H$\textit{\ }$\in $\textit{\ }$L^{\infty }(\Omega
). $\textit{\ Let }$u,$\textit{\ }$v$\textit{\ }$\in $\textit{\ }$%
C^{2}(\Omega ) $ $\cap $ $C^{0}(\bar{\Omega})$ $\cap $ $W^{1,1}(\Omega )$ 
\textit{be two minimizers for }$\mathcal{F}_{H}\mathcal{(\cdot )}$\textit{\
in }$W^{1,1}(\Omega )$\textit{\ (see (\ref{eqn1.0}) and Definition 3.1 in 
\cite{chy}) with }$u$\textit{\ }$=$\textit{\ }$v$\textit{\ on }$\partial
\Omega .$\textit{\ Suppose }$m$\textit{\ }$\geq $\textit{\ }$rank(h_{IJ}(p))$%
\textit{\ }$\geq $\textit{\ }$3$\textit{\ for all }$p$\textit{\ }$\in $%
\textit{\ }$\Omega .$\textit{\ Then }$u$\textit{\ }$=$\textit{\ }$v$\textit{%
\ in }$\Omega .$

\bigskip

A weak version of Theorem A (Corollary A.1, resp.) reads as follows:

\bigskip

\textbf{Theorem A}$^{\prime }.$ \textit{Let }$u,$\textit{\ }$v$\textit{\ }$%
\in $\textit{\ }$W^{2}(\Omega )$ \textit{and }$\vec{F}$\textit{\ }$\in $%
\textit{\ }$W^{1}(\Omega )$\textit{\ where }$\Omega $\textit{\ is a domain
in }$R^{m}.$\textit{\ Suppose for some constant }$C$ $>$ $0,$ $|\nabla u+%
\vec{F}|$ $\geq $ $C,$ $|\nabla v+\vec{F}|$ $\geq $ $C$,\textit{\ }$N_{\vec{F%
}}(u)$\textit{\ }$=$\textit{\ }$N_{\vec{F}}(v)$ \textit{in }$\Omega $ 
\textit{(a.e.)}$.$\textit{\ Assume }%
\begin{equation*}
m\mathit{\ }\geq \mathit{\ }rank(h_{IJ})\mathit{\ }\geq \mathit{\ }3.
\end{equation*}%
\noindent \textit{in\ }$\Omega $ \textit{(a.e.)}$.$ \textit{Then }$\nabla u$%
\textit{\ }$=$\textit{\ }$\nabla v$ \textit{in }$\Omega $ \textit{(a.e.)}$.$

\bigskip

\textbf{Corollary A}$^{\prime }$\textbf{.1}. \textit{Let }$\Omega $\textit{\
be a bounded domain in }$R^{m}.$\textit{\ Take }$\vec{F}$\textit{\ }$\in $%
\textit{\ }$W^{1}(\bar{\Omega})$\textit{\ and }$H$\textit{\ }$\in $\textit{\ 
}$L^{\infty }(\Omega ).$\textit{\ Let }$u,$\textit{\ }$v$\textit{\ }$\in $%
\textit{\ }$W^{2,1}(\Omega )$ \textit{be two minimizers for }$\mathcal{F}_{H}%
\mathcal{(\cdot )}$\textit{\ in }$W^{1,1}(\Omega )$\textit{\ (see (\ref%
{eqn1.0}) and Definition 3.1 in \cite{chy}) with }$u$\textit{\ }$-$\textit{\ 
}$v$\textit{\ }$\in $ $W_{0}^{2,1}(\Omega ).$\textit{\ Suppose for some
constant }$C$ $>$ $0,$ $|\nabla u+\vec{F}|$ $\geq $ $C,$ $|\nabla v+\vec{F}|$
$\geq $ $C$. \textit{Suppose }$m$\textit{\ }$\geq $\textit{\ }$rank(h_{IJ})$%
\textit{\ }$\geq $\textit{\ }$3$\textit{\ in\ }$\Omega $ \textit{(a.e.)}$.$%
\textit{\ Then }$u$\textit{\ }$=$\textit{\ }$v$\textit{\ in }$\Omega $ 
\textit{(a.e.)}$.$

\bigskip

Next we find a nonintegrability condition for $N_{\vec{F}}(u)$ $=$ $N_{\vec{F%
}}(v)$ to imply $\nabla u$ $=$ $\nabla v.$ Let $\Theta _{w}$ $:=$ $%
dw+F_{I}dx_{I}$ for a real function $w$ defined on $\Omega .$ If the
distribution defined by $\Theta _{w}$\ $=$\ $0$\ in $\Omega $ is integrable,
then we have%
\begin{equation*}
\Theta _{w}\wedge d\Theta _{w}=0
\end{equation*}%
\noindent We say $\Theta _{w}$ is integrable (nonintegrable, respectively)
at a point $p$ $\in $ $\Omega $ if $\Theta _{w}\wedge d\Theta _{w}=0$ ($%
\Theta _{w}\wedge d\Theta _{w}\neq 0,$ respectively$)$ at $p.$ The
integrability condition can be described in terms of $h_{IJ}$ and $\nu
_{K}^{w}$ :$=$ ($\partial _{K}w+F_{K})/|\nabla w+\vec{F}|$ (see (\ref{A17});
note that $N_{\vec{F}}(w)$ $=$ $\nu ^{w}$). For $w$ $\in $ $W^{2}(\Omega )$%
\textit{\ }and\textit{\ }$\vec{F}\mathit{\ }\in \mathit{\ }W^{1}(\Omega ),$
we say $\Theta _{w}$ is nonintegrable if $\Theta _{w}\wedge d\Theta _{w}\neq
0$ in $\Omega $ a.e..

\bigskip

\textbf{Theorem B}. \textit{Let }$u,$\textit{\ }$v$\textit{\ }$\in $\textit{%
\ }$C^{2}(\Omega )$ \textit{and }$\vec{F}$\textit{\ }$\in $\textit{\ }$%
C^{1}(\Omega )$\textit{\ where }$\Omega $\textit{\ is a domain in }$R^{m}.$%
\textit{\ Suppose }$N_{\vec{F}}(u)$\textit{\ }$=$\textit{\ }$N_{\vec{F}}(v)$%
\textit{\ in }$\Omega $\textit{\ }$\backslash $\textit{\ }$[S_{\vec{F}}(u)$%
\textit{\ }$\cup $\textit{\ }$S_{\vec{F}}(v)].$\textit{\ Suppose for each
point in }$\Omega $\textit{\ }$\backslash $\textit{\ }$[S_{\vec{F}}(u)$%
\textit{\ }$\cup $\textit{\ }$S_{\vec{F}}(v)],$ \textit{either }$\Theta _{u}$%
\textit{\ is nonintegrable or }$\Theta _{v}$\textit{\ is nonintegrable}$.$%
\textit{\ Then }$\nabla u$\textit{\ }$=$\textit{\ }$\nabla v$ \textit{in} $%
\Omega $\textit{\ }$\backslash $\textit{\ }$[S_{\vec{F}}(u)$\textit{\ }$\cup 
$\textit{\ }$S_{\vec{F}}(v)].$

\bigskip

Again by adding the boundary condition we then have the uniqueness of
minimizers for $\mathcal{F}_{H}.$

\bigskip

\textbf{Corollary B.1}. \textit{Let }$\Omega $\textit{\ be a bounded domain
in }$R^{m}.$\textit{\ Take }$\vec{F}$\textit{\ }$\in $\textit{\ }$C^{1}(\bar{%
\Omega})$\textit{\ and }$H$\textit{\ }$\in $\textit{\ }$L^{\infty }(\Omega
). $\textit{\ Let }$u,$\textit{\ }$v$\textit{\ }$\in $\textit{\ }$%
C^{2}(\Omega ) $ $\cap $ $C^{0}(\bar{\Omega})$ $\cap $ $W^{1,1}(\Omega )$ 
\textit{be two minimizers for }$\mathcal{F}_{H}\mathcal{(\cdot )}$\textit{\
in }$W^{1,1}(\Omega )$\textit{\ (see Definition 3.1 in \cite{chy}) with }$u$%
\textit{\ }$=$\textit{\ }$v$\textit{\ on }$\partial \Omega .$\textit{\
Suppose for each point in }$\Omega $\textit{\ }$\backslash $\textit{\ }$[S_{%
\vec{F}}(u)$\textit{\ }$\cup $\textit{\ }$S_{\vec{F}}(v)],$ \textit{either }$%
\Theta _{u}$\textit{\ is nonintegrable or }$\Theta _{v}$\textit{\ is
nonintegrable}$.$\textit{\ Then }$u$\textit{\ }$=$\textit{\ }$v$\textit{\ in 
}$\Omega .$

\bigskip

A weak version of Theorem B (Corollary B.1, resp.) reads as follows:

\bigskip

\textbf{Theorem B}$^{\prime }$. \textit{Let }$u,$\textit{\ }$v$\textit{\ }$%
\in $\textit{\ }$W^{2}(\Omega )$ \textit{and }$\vec{F}$\textit{\ }$\in $%
\textit{\ }$W^{1}(\Omega )$\textit{\ where }$\Omega $\textit{\ is a domain
in }$R^{m}.$\textit{\ Suppose for some constant }$C$ $>$ $0,$ $|\nabla u+%
\vec{F}|$ $\geq $ $C,$ $|\nabla v+\vec{F}|$ $\geq $ $C$, $N_{\vec{F}}(u)$%
\textit{\ }$=$\textit{\ }$N_{\vec{F}}(v)$\textit{\ in }$\Omega $\textit{\
(a.e.)}$.$\textit{\ Suppose} \textit{either }$\Theta _{u}$\textit{\ is
nonintegrable or }$\Theta _{v}$\textit{\ is nonintegrable}$.$\textit{\ Then }%
$\nabla u$\textit{\ }$=$\textit{\ }$\nabla v$ \textit{in} $\Omega $\textit{\
(a.e.)}$.$

\bigskip

\textbf{Corollary B}$^{\prime }$\textbf{.1}. \textit{Let }$\Omega $\textit{\
be a bounded domain in }$R^{m}.$\textit{\ Take }$\vec{F}$\textit{\ }$\in $%
\textit{\ }$W^{1}(\Omega )$\textit{\ and }$H$\textit{\ }$\in $\textit{\ }$%
L^{\infty }(\Omega ).$\textit{\ Let }$u,$\textit{\ }$v$\textit{\ }$\in $%
\textit{\ }$W^{2,1}(\Omega )$ \textit{be two minimizers for }$\mathcal{F}_{H}%
\mathcal{(\cdot )}$\textit{\ in }$W^{1,1}(\Omega )$\textit{\ (see Definition
3.1 in \cite{chy}) with }$u$\textit{\ }$-$\textit{\ }$v$\textit{\ }$\in $ $%
W_{0}^{2,1}(\Omega ).$\textit{\ Suppose for some constant }$C$ $>$ $0,$ $%
|\nabla u+\vec{F}|$ $\geq $ $C,$ $|\nabla v+\vec{F}|$ $\geq $ $C$ \textit{%
(a.e.).} \textit{Suppose either }$\Theta _{u}$\textit{\ is nonintegrable or }%
$\Theta _{v}$\textit{\ is nonintegrable}$.$\textit{\ Then }$u$\textit{\ }$=$%
\textit{\ }$v$\textit{\ in }$\Omega $ \textit{(a.e.)}$.$

\bigskip

Note that in the above results the dimension $"m"$ is not necessarily even.
We can also extend Theorem 1.1 under a condition more general than $\func{div%
}\vec{F}^{\ast }$ $>$ (or $<)$ $0$ while the dimension $"m"$ is not
necessarily even. Define $\vec{G}^{b}$ for $\vec{G}$ $=$ $(G_{1},$ $...,$ $%
G_{m})$ by

\begin{equation*}
\vec{G}^{b}:=(\sum_{k=1}^{m}a^{1k}G_{k},\sum_{k=1}^{m}a^{2k}G_{k},...,%
\sum_{k=1}^{m}a^{mk}G_{k})
\end{equation*}

\noindent where $a^{jk\prime }s$ are real constants such that $a^{jk}+a^{kj}$
$=$ $0$ for $1$ $\leq $ $j,k$ $\leq $ $m.$ Note that $\vec{G}^{b}$ $=$ $\vec{%
G}^{\ast }$ for $m$ $=$ $2n,$ $a^{2j-1,2j}$ $=$ $-a^{2j,2j-1}$ $=$ $1,$ $%
1\leq j\leq n,$ $a^{jk}$ $=$ $0$ otherwise.

\bigskip

\textbf{Theorem C}. \textit{Let} $\Omega $\textit{\ be a bounded domain in }$%
R^{m}.$ \textit{Let }$u,$\textit{\ }$v$\textit{\ }$\in $\textit{\ }$%
W^{1,2}(\Omega )$\textit{\ be two minimizers for }$\mathcal{F}_{H}$\textit{\
such that }$u-v$\textit{\ }$\in $\textit{\ }$W_{0}^{1,2}(\Omega ).$\textit{\
Suppose }$H\in L^{\infty }(\Omega )$\textit{\ and }$\vec{F}\in
W^{1,2}(\Omega )$\textit{\ satisfying }%
\begin{equation*}
\func{div}\vec{F}^{b}=\sum_{j,k=1}^{m}a^{jk}\partial _{j}F_{k}>0\text{ (}<0,%
\text{ resp.) (a.e.)}
\end{equation*}%
\textit{\noindent where }$a^{jk}$\textit{'s are real constants such that }$%
a^{jk}+a^{kj}=0.$\textit{\ Then }$u\equiv v$\textit{\ in }$\Omega $\textit{\
(a.e.).}

\bigskip

Compare (\ref{A8}) with the Euclidean situation: $\delta _{I}\nu _{J}-\delta
_{J}\nu _{I}=0$ where 
\begin{equation*}
\nu =(\nu _{I})=\frac{(-\nabla u,1)}{\sqrt{1+|\nabla u|^{2}}}
\end{equation*}

\noindent denotes the unit normal to the graph defined by $u.$ It is a known
fact that $\nu $ can be realized as the unit normal vectors of a family of
(hyper)surfaces filling up a region if and only if $\delta _{I}\nu
_{J}-\delta _{J}\nu _{I}$ $=$ $0$ (see page 3 in \cite{MM})$.$

Recall that in our situation, the horizontal normal $\nu ^{u}$ of a graph
defined by $u$ reads%
\begin{equation*}
\nu ^{u}=\frac{\nabla u+\vec{F}}{|\nabla u+\vec{F}|}.
\end{equation*}

\noindent Let $\nu _{K}^{u}$ denote the components of $\nu ^{u}.$ Recall
that the horizontal tangential operator $\delta _{K}^{u}$ is defined by $%
\delta _{K}^{u}$ $=$ $\partial _{K}-\nu _{K}^{u}\nu _{J}^{u}\partial _{J}$
(cf. (\ref{A2}))$.$ Recall that $D_{u}$ $:=$ $|\nabla u+\vec{F}|$ and $%
h_{IJ} $ $:=$ $\partial _{I}F_{J}-\partial _{J}F_{I}$ (cf. (\ref{A7}))$.$ In
Section 2, we deduce%
\begin{eqnarray*}
&&\delta _{I}^{u}\nu _{J}^{u}-\delta _{J}^{u}\nu _{I}^{u} \\
&=&\frac{1}{D_{u}}\{h_{IJ}-\nu _{J}^{u}\nu _{K}^{u}h_{IK}-\nu _{I}^{u}\nu
_{K}^{u}h_{KJ}\}.
\end{eqnarray*}

\noindent (i.e., (\ref{A8})). In view of the Euclidean situation, we ask the
following question:

\bigskip

\textbf{Question D}: Given $D>0,$ a unit vector $\nu =(\nu _{J}),$ $\vec{F}%
=(F_{J})$ locally in $R^{m}$ satisfying%
\begin{eqnarray}
&&\delta _{I}\nu _{J}-\delta _{J}\nu _{I}  \label{1.1.4} \\
&=&\frac{1}{D}\{h_{IJ}-\nu _{J}\nu _{K}h_{IK}-\nu _{I}\nu _{K}h_{KJ}\}, 
\notag
\end{eqnarray}

\noindent can we find $u$ such that $\nu $ $=$ $\frac{\nabla u+\vec{F}}{D}?$

\bigskip

Here $\delta _{J}$ :$=$ $\partial _{J}-\nu _{J}\nu _{K}\partial _{K}$. In
fact, we are asking if (\ref{1.1.4}) is an integrability condition for $\nu $
to be the horizontal normal of a graph defined by $u.$ It turns out that we
need a condition other than (\ref{1.1.4}) to conclude $\nu $ $=$ $\frac{%
\nabla u+\vec{F}}{D}$. Denote 1-forms $\nu _{J}dx^{J},$ $F_{J}dx^{J}$ by $%
\nu ^{\#},$ $F^{\#},$ resp.. When $\nu $ $=$ $\frac{\nabla u+\vec{F}}{D},$
we get $D\nu ^{\#}$ $=$ $du+F^{\#}.$ It follows that $d(D\nu ^{\#})$ $=$ $%
dF^{\#}.$ So we have%
\begin{equation}
\nu ^{\#}\lrcorner d(D\nu ^{\#}-F^{\#})=0.  \label{1.1.5}
\end{equation}

Recall that the interior product $\eta \lrcorner \omega $ of 1-form $\eta $
and 2-form $\omega $ is defined to be ($\eta _{\#})\lrcorner \omega $ $:=$ $%
\omega (\eta _{\#})$ where $\eta _{\#}$ is the corresponding vector of $\eta 
$ with respect to the Euclidean metric. In practice, we have $dx^{i}$ $%
\lrcorner $ $(dx^{j}\wedge dx^{k})$ $=$ $<dx^{i},dx^{j}>$ $dx^{k}$ $-$ $%
<dx^{i},dx^{k}>$ $dx^{j}$ $=$ $\delta ^{ij}dx^{k}$ $-$ $\delta ^{ik}dx^{j}$ (%
$\delta ^{ij}$ denotes the Kronecker delta, i.e., $\delta ^{ij}$ $=$ $1$ if $%
i$ $=$ $j;$ $\delta ^{ij}\ =$ $0$ if $i$ $\neq $ $j),$ and hence ($\eta
_{i}dx^{i})$ $\lrcorner $ $(\omega _{jk}dx^{j}\wedge dx^{k})$ $=$ $\eta
_{i}\omega _{jk}\{dx^{i}$ $\lrcorner $ $(dx^{j}\wedge dx^{k})\}$ $=$ $\eta
_{j}\omega _{jk}dx^{k}$ $-$ $\eta _{k}\omega _{jk}dx^{j}$ $=$ $\eta
_{j}(\omega _{jk}$ $-$ $\omega _{kj})dx^{k}.$

We can view (\ref{1.1.5}) as a system of first order equations in $D$
coupled with (\ref{1.1.4}), a system of first order equations in $\nu .$ It
is not hard to rewrite (\ref{1.1.5}) as follows:%
\begin{equation}
\delta _{K}D=\nu _{J}(\partial _{J}\nu _{K})D-\nu _{J}h_{JK}  \label{1.1.6}
\end{equation}

\noindent for any $K$ (see (\ref{D6}) in Section 4). From the above
discussion we learn that (\ref{1.1.4}) and (\ref{1.1.5}) (or equivalently, (%
\ref{1.1.6})) are two necessary conditions for $\nu $ to be the horizontal
normal associated to a function $u,$ i.e., $\nu $ $=$ $\frac{\nabla u+\vec{F}%
}{D}.$ Conversely, they are also sufficient as we answer Question D in the
following integrability theorem. For simplicity, we work in $C^{\infty }$
category for this problem.

\bigskip

\textbf{Theorem E}. \textit{Given (}$C^{\infty }$\textit{\ smooth) }$D>0,$%
\textit{\ a unit vector }$\nu =(\nu _{J}),$\textit{\ }$\vec{F}=(F_{J})$%
\textit{\ locally in }$R^{m},$\textit{\ }$m$ $\geq $ $2,$ \textit{satisfying
(\ref{1.1.4}), and (\ref{1.1.5}) or equivalently (\ref{1.1.6}), we can find
a (}$C^{\infty }$\textit{\ smooth) function }$u$\textit{\ locally such that }%
$\nu $\textit{\ }$=$\textit{\ }$\frac{\nabla u+\vec{F}}{D}.$

\bigskip

If $m$ $=$ $2,$ (\ref{1.1.5}) or (\ref{1.1.6}) is equivalent to $d(D\nu
^{\#}-F^{\#})$ $=$ $0$ (see Proposition 4.1)$.$ However, for higher
dimensions, we can have the situation that (\ref{1.1.5}) holds while $d(D\nu
^{\#}-F^{\#})$ $\neq $ $0$ (see Example 4.2). Also we can have the situation
that (\ref{1.1.4}) holds while (\ref{1.1.5}) does not hold (see Example 4.3).

Comparing with the fundamental theorem for surfaces in the 3-dimensional
Heisenberg group in \cite{chmy2}, we don't prescribe $p$-mean curvature $H$
here, but prescribe arbitrary $\vec{F}$ instead of fixed $\vec{F}$ $=$ $(-y,$
$x)$ in \cite{chmy2}. Equation (\ref{1.1.5}) or (\ref{1.1.6}) corresponds to
a Codazzi-like equation (cf. (1.17) in \cite{chmy2}). See (\ref{D11.1}) and
the discussion before Example 4.2 in Section 4 (see also a recent preprint
of Hung-Lin Chiu \cite{Chiu}).

The idea of proof for Theorem E is to show that $U_{I}$ $:=$ $D\nu
_{I}-F_{I} $ satisfy the integrability condition $\partial _{I}U_{J}$ $=$ $%
\partial _{J}U_{I}$ (and hence $U_{I}$ $=$ $\partial _{I}u$ for some
function $u).$ Let $U_{IJ}$ $:=$ $\partial _{I}U_{J}$ $-$ $\partial
_{J}U_{I}.$ A direct computation shows that%
\begin{equation}
U_{IJ}-\nu _{J}\nu _{K}U_{IK}-\nu _{I}\nu _{K}U_{KJ}=0  \label{1.1.7}
\end{equation}

\noindent due to condition (\ref{1.1.4}). Observe that in terms of
differential forms, we can write (\ref{1.1.7}) as follows:%
\begin{eqnarray}
d(U_{I}dx^{I}) &=&\frac{1}{2}U_{IJ}dx^{I}\wedge dx^{J}  \label{1.1.8} \\
&=&(U\nu )^{\#}\wedge \nu ^{\#}  \notag
\end{eqnarray}

\noindent where $U$ denotes the matrix $(U_{IJ})$ and $\nu $ is viewed as a
column vector in $U\nu .$ We then observe that $D\nu ^{\#}$ $-$ $F^{\#}$ $=$ 
$U_{I}dx^{I},$ and hence $d(D\nu ^{\#}-F^{\#})$ $=$ $\frac{1}{2}%
U_{IJ}dx^{I}\wedge dx^{J}.$ Substituting (\ref{1.1.8}) into condition (\ref%
{1.1.5}) (or (\ref{1.1.6})) gives ($U\nu )^{\#}$ $=$ $0.$ By (\ref{1.1.8})
again, we get $U_{IJ}$ $=$ $0,$ i.e., $\partial _{I}U_{J}$ $=$ $\partial
_{J}U_{I}.$

\bigskip

\section{\textbf{Proofs of Theorems A and B}}

Recall in \cite{chy} that $S_{\vec{F}}(u)$ $:=$ $\{p$ $\in $ $\Omega $ 
\TEXTsymbol{\vert} ($\nabla u+\vec{F})(p)$ $=$ $0\}$ and $N_{\vec{F}}(u)$ :$%
= $ $\frac{\nabla u+\vec{F}}{|\nabla u+\vec{F}|}.$ The idea of the proof for
the uniqueness in \cite{chy} is to show that $N_{\vec{F}}(u)$ $=$ $N_{\vec{F}%
}(v)$ in $\Omega $ $\backslash $ $[S_{\vec{F}}(u)$ $\cup $ $S_{\vec{F}}(v)]$
(a.e.) first for two minimizers $u,$ $v,$ say$,$ in $W^{1,1}(\Omega )$ such
that $u-v$ $\in $ $W_{0}^{1.1}(\Omega ).$ Then to show that $\nabla u$ $=$ $%
\nabla v$ (and hence $u$ $=$ $v),$ we invoke an equality (see (5.3) in \cite%
{chy}) and an argument of integrating by parts (see Theorem 5.3 in \cite{chy}%
). To make this approach work, we need to assume $m$ $=$ $2n$ and $\func{div}%
\vec{F}^{\ast }$ $>$ (or $<,$ resp.) $0$ (a.e.). In this section, we are
going to give another approach to show that $N_{\vec{F}}(u)$ $=$ $N_{\vec{F}%
}(v)$ implies $\nabla u$ $=$ $\nabla v.$ Note that in this approach, we do
not need to assume $m$ $=$ $2n.$

To explain the idea, we assume $u,$ $v$ $\in $ $C^{2}.$ Write $\nabla u$ $=$ 
$(u_{K})$ where $u_{K}$ $:=$ $\partial _{K}u,$ $\partial _{K}$ $:=$ $\frac{%
\partial }{\partial x_{K}},$ $K$ $=$ $1,$ $2,$ $...,$ $m$ and $\vec{F}$ $=$ $%
(F_{K}).$ So we can write%
\begin{equation*}
N_{\vec{F}}(u)=(\frac{u_{K}+F_{K}}{D_{u}})
\end{equation*}

\noindent where $D_{u}$ $:=$ $|\nabla u+\vec{F}|.$ Let $\hat{u}_{K}$ $:=$ $%
u_{K}+F_{K}$ and $\nu _{K}^{u}$ $:=$ $\frac{\hat{u}_{K}}{D_{u}},$ components
of $N(u).$ Define the horizontal tangential operator $\delta _{K}^{u}$ by%
\begin{equation}
\delta _{K}^{u}=\partial _{K}-\nu _{K}^{u}\nu _{J}^{u}\partial _{J}
\label{A2}
\end{equation}

\noindent (summing over $J$; summation convention hereafter$).$ We compute%
\begin{equation}
\partial _{I}\nu _{J}^{u}=\frac{\partial _{I}\hat{u}_{J}}{D_{u}}-\frac{\hat{u%
}_{J}\hat{u}_{K}\partial _{I}\hat{u}_{K}}{D_{u}^{3}}.  \label{A3}
\end{equation}

\noindent Hence from (\ref{A2}), (\ref{A3}), and the definition of $\nu
_{J}^{u},$ we have%
\begin{eqnarray}
\delta _{I}^{u}\nu _{J}^{u} &=&\partial _{I}\nu _{J}^{u}-\nu _{I}^{u}\nu
_{K}^{u}\partial _{K}\nu _{J}^{u}  \label{A4} \\
&=&\frac{\partial _{I}\hat{u}_{J}}{D_{u}}-\frac{\hat{u}_{J}\hat{u}%
_{K}\partial _{I}\hat{u}_{K}}{D_{u}^{3}}-\frac{\hat{u}_{I}\hat{u}%
_{K}\partial _{K}\hat{u}_{J}}{D_{u}^{3}}  \notag \\
&&+\frac{\hat{u}_{I}\hat{u}_{K}\hat{u}_{J}\hat{u}_{L}\partial _{K}\hat{u}_{L}%
}{D_{u}^{5}}.  \notag
\end{eqnarray}

We can now compute%
\begin{eqnarray}
&&\delta _{I}^{u}\nu _{J}^{u}-\delta _{J}^{u}\nu _{I}^{u}  \label{A5} \\
&=&\frac{\partial _{I}\hat{u}_{J}-\partial _{J}\hat{u}_{I}}{D_{u}}-\frac{%
\hat{u}_{J}\hat{u}_{K}(\partial _{I}\hat{u}_{K}-\partial _{K}\hat{u}_{I})}{%
D_{u}^{3}}  \notag \\
&&-\frac{\hat{u}_{I}\hat{u}_{K}(\partial _{K}\hat{u}_{J}-\partial _{J}\hat{u}%
_{K})}{D_{u}^{3}}  \notag
\end{eqnarray}

\noindent by (\ref{A4}) and noting that the term involving $D_{u}^{5}$ is
symmetric in $I,$ $J.$ From the definition of $\hat{u}_{J}$ we have%
\begin{eqnarray}
&&\partial _{I}\hat{u}_{J}-\partial _{J}\hat{u}_{I}  \label{A6} \\
&=&\partial _{I}(u_{J}+F_{J})-\partial _{J}(u_{I}+F_{I})  \notag \\
&=&\partial _{I}F_{J}-\partial _{J}F_{I}.  \notag
\end{eqnarray}

\noindent Let%
\begin{equation}
h_{IJ}:=\partial _{I}F_{J}-\partial _{J}F_{I}.  \label{A7}
\end{equation}

\noindent So in view of (\ref{A6}) and (\ref{A7}), we can write (\ref{A5})
as follows:%
\begin{eqnarray}
&&\delta _{I}^{u}\nu _{J}^{u}-\delta _{J}^{u}\nu _{I}^{u}  \label{A8} \\
&=&\frac{1}{D_{u}}\{h_{IJ}-\nu _{J}^{u}\nu _{K}^{u}h_{IK}-\nu _{I}^{u}\nu
_{K}^{u}h_{KJ}\}.  \notag
\end{eqnarray}

Now suppose $N_{\vec{F}}(u)$ $=$ $N_{\vec{F}}(v).$ Then $\nu _{K}^{u}$ $=$ $%
\nu _{K}^{v}$ and $\delta _{K}^{u}$ $=$ $\delta _{K}^{v}$ from the
definition. It follows from (\ref{A8}) that%
\begin{equation}
\{h_{IJ}-\nu _{J}^{u}\nu _{K}^{u}h_{IK}-\nu _{I}^{u}\nu _{K}^{u}h_{KJ}\}(%
\frac{1}{D_{u}}-\frac{1}{D_{v}})=0  \label{A9}
\end{equation}

\noindent (in $\Omega $ $\backslash $ $[S_{\vec{F}}(u)$ $\cup $ $S_{\vec{F}%
}(v)]).$ If $\nabla u$ $\neq $ $\nabla v,$ then $D_{u}$ $\neq $ $D_{v}.$
Therefore we have%
\begin{equation}
h_{IJ}-\nu _{J}^{u}\nu _{K}^{u}h_{IK}-\nu _{I}^{u}\nu _{K}^{u}h_{KJ}=0
\label{A10}
\end{equation}

\noindent by (\ref{A9}). Observe that $h$ $=$ $(h_{IJ})$ is a skew-symmetric
matrix by (\ref{A7}), i.e.%
\begin{equation}
h+h^{T}=0\text{ or }h_{IJ}+h_{JI}=0  \label{A11}
\end{equation}

\noindent where $h^{T}$ denotes the transpose of $h.$

\bigskip

\textbf{Lemma 2.1.} \textit{Suppose }$h$\textit{\ is a skew-symmetric real }$%
m\times m$\textit{\ matrix (}$m$\textit{\ }$\geq $\textit{\ }$2)$\textit{\
such that}%
\begin{equation}
h-h\nu \nu ^{T}-\nu \nu ^{T}h=0  \label{A12}
\end{equation}

\noindent \textit{where }$\nu $\textit{\ is a (}$m\times 1)$\textit{\ unit
column real vector and }$\nu ^{T}$\textit{\ is the transpose of }$\nu ,$%
\textit{\ a (}$1\times m)$\textit{\ unit row vector. Then we have}

\begin{equation}
rank(h)=0\text{ or }2.  \label{A12.1}
\end{equation}

\noindent \textit{where }$rank(h)$\textit{\ denotes the rank of }$h.$

\bigskip

\proof
Multiplying (\ref{A12}) by $h$ and then taking the trace, we obtain%
\begin{eqnarray}
Tr(h^{2}) &=&Tr(h\nu \nu ^{T}h)+Tr(\nu \nu ^{T}hh)  \label{A13} \\
&=&2Tr(h\nu \nu ^{T}h)\text{ (since }Tr(\nu \nu ^{T}hh)=Tr(h\nu \nu ^{T}h)) 
\notag \\
&=&-2||h\nu ||^{2}\text{ (since }h^{T}=-h\text{ by (\ref{A11})}).  \notag
\end{eqnarray}

\noindent Here $||$ $\cdot $ $||$ denotes the Euclidean norm. Observe that
the eigenvalues of $h$ (being skew-symmetric) are purely imaginary and if $%
i\lambda $ ($\lambda $ $\in $ $R\backslash \{0\})$ is a nonzero eigenvalue
(with an eigenvector $w)$, then $-i\lambda $ is also an eigenvalue (with an
eigenvector independent of $w).$ It follows that $h^{2}$ has an eigenvalue $%
-\lambda ^{2}$ of multiplicity $2.$ Let $i\lambda _{1},$ $i\lambda _{2},$
..,.$i\lambda _{k}$ ($\lambda _{j}$ $\in $ $R\backslash \{0\}),$ $|\lambda
_{1}|$ $\geq $ $|\lambda _{2}|$ $\geq $ ...$\geq $ $|\lambda _{k}|$ $>$ $0,$ 
$2k$ $\leq $ $m,$ be all nonzero eigenvalues of $h$ (if $h$ $\neq $ $0$)
while $-\lambda _{j}^{2},$ $j$ $=$ $1,$ $2,$ ...,$k,$ are all nonzero
eigenvalues of $h^{2}$ (each of which has multiplicity 2)$.$ From (\ref{A13}%
) we can easily get

\begin{equation*}
2\lambda _{1}^{2}\leq 2\sum_{j=1}^{k}\lambda _{j}^{2}=2||h\nu ||^{2}\leq
2\lambda _{1}^{2}.
\end{equation*}

\noindent Therefore $\ k$ $=$ $1$ and hence ($\ref{A12.1})$ follows.

\endproof%

\bigskip

\proof
\textbf{(of Theorem A)} This follows from the above discussion. Take $h$ $=$ 
$(h_{IJ}(p))$ where $h_{IJ}$ $=$ $\partial _{I}F_{J}-\partial _{J}F_{I}$ and 
$\nu $ $=$ ($\nu _{J}^{u}(p))$ in Lemma 2.1. There holds (\ref{A12}) if $%
\nabla u$ $\neq $ $\nabla v$ at $p$ by (\ref{A9}). Now the conclusion of
Lemma 2.1 contradicts the assumption on $rank(h).$

\endproof%

\bigskip

We remark that $N_{\vec{F}}(u)$\textit{\ }$=$\textit{\ }$N_{\vec{F}}(v)$ (in
a region) does not imply $\nabla u$\textit{\ }$=$\textit{\ }$\nabla v$ in
dimension 2 as shown in the following example.

\bigskip

\textbf{Example 2.2.} Take $u$ $=$ $xy$ and $v$ $=$ $xy$ $+$ $y$ which
define graphs over the $xy$-plane in the Heisenberg group of dimension 3. So
in this situation we have $\vec{F}$ $=$ $(-y,x)$. It is straightforward to
compute $\nabla u$ $=$ $(y,x),$ $\nabla v$ $=$ $(y,x+1),$ $\nabla u+\vec{F}$ 
$=$ $(0,2x),$ $\nabla v+\vec{F}$ $=$ $(0,2x+1).$ We observe that $N_{\vec{F}%
}(u)$\textit{\ }$=$\textit{\ }$N_{\vec{F}}(v)$ $=$ $(0,1)$ in the region $%
\{x>0\}$. On the other hand, it is clear that $\nabla u$ $\neq $ $\nabla v.$

\bigskip

\proof
\textbf{(of Corollary A.1)} By Theorem 5.1 in \cite{chy}, we have $N_{\vec{F}%
}(u)$\textit{\ }$=$\textit{\ }$N_{\vec{F}}(v)$ in $\Omega $\textit{\ }$%
\backslash $\textit{\ }$[S_{\vec{F}}(u)$\textit{\ }$\cup $\textit{\ }$S_{%
\vec{F}}(v)]$ (taking regular values $\varepsilon _{j}\rightarrow 0$ and $1$%
, resp.). Then it follows from Theorem A that%
\begin{equation}
\nabla u=\nabla v\text{ in}\mathit{\ }\Omega \mathit{\ }\backslash \mathit{\ 
}[S_{\vec{F}}(u)\mathit{\ }\cup \mathit{\ }S_{\vec{F}}(v)].  \label{A13.1}
\end{equation}%
\noindent We claim that both $S_{\vec{F}}(u)$ and $S_{\vec{F}}(v)$ are
nowhere dense in $\Omega .$ Suppose the converse holds. Then we can find a
small ball $B$ contained in either $S_{\vec{F}}(u)$ or $S_{\vec{F}}(v),$ say 
$B$ $\subset $ $S_{\vec{F}}(u).$ This means that $\nabla u$ $+$ $\vec{F}$ $=$
$0$ in $B.$ It follows that $F_{I}$ $=$ $-u_{I}$ and hence $\partial
_{I}F_{J}$ $=$ $-\partial _{I}u_{J}$ $=$ $-\partial _{J}u_{I}$ $=$ $\partial
_{J}F_{I}.$ Therefore $h_{IJ}$ $=$ $\partial _{I}F_{J}-\partial _{J}F_{I}$ $%
= $ $0$ in $B$ for all $I,J,$ contradicting the condition on the rank of $%
(h_{IJ}).$ The above argument also works for $B$ $\subset $ $S_{\vec{F}}(v).$
So we have shown that both $S_{\vec{F}}(u)$ and $S_{\vec{F}}(v)$ are nowhere
dense in $\Omega .$ It follows that $S_{\vec{F}}(u)$\textit{\ }$\cup $%
\textit{\ }$S_{\vec{F}}(v)$ is nowhere dense in $\Omega $. Therefore by (\ref%
{A13.1}) we have $u$ $-$ $v$ $=$ $c,$ a constant in $\Omega .$ Since $u$ and 
$v$\textit{\ }are continous up to the boundary $\partial \Omega $ and $u$ $=$%
\textit{\ }$v$\textit{\ }on\textit{\ }$\partial \Omega ,$ we have $c$ $=$ $%
0. $

\endproof%

\bigskip

We remark that for $m$ $=$ $2,$ $\vec{F}$\textit{\ }$\in $\textit{\ }$%
C^{1}(\Omega )$ and $w$ $\in $ $C^{1}(\Omega )$, $S_{\vec{F}}(w)$ is nowhere
dense in $\Omega $ if $\func{div}\vec{F}^{\ast }$ $>$ $0$ (or $<$ $0,$ resp.$%
)$ in $\Omega $ (cf. Lemma 3.1 in \cite{chmy2}; in fact, we can extend this
result to $m$ $=$ $2n$ or even general dimensions, see Proposition 3.4\ in
this paper)$.$ Also note that the size of the singular set can be measured
in terms of the rank of ($h_{IJ})$ (see Theorem D in \cite{chy} where we
need to assume $u$ $\in $ $C^{2},$ $\vec{F}$\textit{\ }$\in $\textit{\ }$%
C^{1}$ in view of Balogh's $C^{1,1}$ examples in \cite{Ba}).

We can interpret (\ref{A10}) as an integrability condition for hypersurfaces
annihilated by the one-form%
\begin{eqnarray}
\Theta _{u} &:&=du+F_{I}dx^{I}  \label{A14} \\
&=&(u_{I}+F_{I})dx^{I}  \notag \\
&=&D_{u}\nu _{I}^{u}dx^{I}.  \notag
\end{eqnarray}

\bigskip

\textbf{Lemma 2.3}. \textit{Let }$u$\textit{\ }$\in $\textit{\ }$%
C^{2}(\Omega )$\textit{\ and }$\vec{F}$\textit{\ }$\in $\textit{\ }$%
C^{1}(\Omega )$\textit{\ where }$\Omega $\textit{\ is a domain of }$R^{m}.$ 
\textit{Then} \textit{in }$\Omega $\textit{\ }$\backslash $\textit{\ }$S_{%
\vec{F}}(u),$\textit{\ }$\Theta _{u}$\textit{\ is integrable (meaning the
distribution defined by }$\Theta _{u}$\textit{\ }$=$\textit{\ }$0$\textit{\
is integrable) if and only if (\ref{A10}) holds.}

\bigskip

\proof
Observe that $\Theta _{u}$ is integrable if and only if $\Theta _{u}\wedge
d\Theta _{u}$ $=$ $0$ by Frobenius' integrability theorem. We then compute%
\begin{eqnarray}
d\Theta _{u} &=&d(du+F_{I}dx^{I})  \label{A15} \\
&=&\partial _{J}F_{I}dx^{J}\wedge dx^{I}  \notag \\
&=&\frac{1}{2}h_{IJ}dx^{I}\wedge dx^{J}  \notag
\end{eqnarray}

\noindent (recall that $h_{IJ}:=\partial _{I}F_{J}-\partial _{J}F_{I})$ and%
\begin{eqnarray}
&&\Theta _{u}\wedge d\Theta _{u}  \label{A16} \\
&=&\frac{1}{2}D_{u}\nu _{K}^{u}h_{IJ}dx^{K}\wedge dx^{I}\wedge dx^{J}  \notag
\end{eqnarray}

\noindent by (\ref{A14}) and (\ref{A15}). So $\Theta _{u}\wedge d\Theta _{u}$
$=$ $0$ if and only if 
\begin{equation}
\nu _{K}^{u}h_{IJ}+\nu _{I}^{u}h_{JK}+\nu _{J}^{u}h_{KI}=0  \label{A17}
\end{equation}

\noindent by (\ref{A16}) and noting that $h_{IJ}$ $=$ $-h_{JI}.$ Multiplying
(\ref{A17}) by $\nu _{K}^{u}$ and summing over $K,$ we obtain (\ref{A10}).
Conversely if (\ref{A10}) holds, we can also easily deduce (\ref{A17}).

\endproof%

\bigskip

\proof
\textbf{(of Theorem B)}

From $N_{\vec{F}}(u)$\textit{\ }$=$\textit{\ }$N_{\vec{F}}(v)$ in $\Omega $%
\textit{\ }$\backslash $\textit{\ }$[S_{\vec{F}}(u)$\textit{\ }$\cup $%
\textit{\ }$S_{\vec{F}}(v)]$ (comparing with the proof of Theorem A), we
know that $\nabla u$ $+$ $\vec{F}$ is parallel to $\nabla v$ $+$ $\vec{F}.$
It follows that in $\Omega $\textit{\ }$\backslash $\textit{\ }$[S_{\vec{F}%
}(u)$\textit{\ }$\cup $\textit{\ }$S_{\vec{F}}(v)],$ there holds%
\begin{equation}
du+F_{I}dx^{I}=\lambda (dv+F_{I}dx^{I})  \label{A18}
\end{equation}

\noindent for some nonzero function $\lambda .$ Subtracting $dv+F_{I}dx^{I}$
from (\ref{A18}) gives%
\begin{equation}
d(u-v)=(\lambda -1)(dv+F_{I}dx^{I}).  \label{A19}
\end{equation}

\noindent Taking exterior differentiation of (\ref{A19}), we obtain%
\begin{equation}
0=d\lambda \wedge \Theta _{v}+(\lambda -1)d\Theta _{v}  \label{A20}
\end{equation}

\noindent where $\Theta _{v}$ $:=$ $dv+F_{I}dx_{I}.$ Wedging (\ref{A20})
with $\Theta _{v}$ we get%
\begin{equation}
(\lambda -1)\Theta _{v}\wedge d\Theta _{v}=0.  \label{A21}
\end{equation}

\noindent Observe that $\lambda $ $=$ $1$ if and only if $\nabla u$ $=$ $%
\nabla v.$ So if $\nabla u$ $\neq $ $\nabla v,$ we have $\Theta _{v}\wedge
d\Theta _{v}$ $=$ $0$ (and $\Theta _{u}\wedge d\Theta _{u}$ $=$ $0,$ resp.)
by (\ref{A21}) (an identity replacing $v$ by $u,$ resp$.)$. This contradicts
the nonintegrability of $\Theta _{v}$ or $\Theta _{u},$ the main assumption
of Theorem B. Therefore we have $\nabla u$ $=$ $\nabla v$ in $\Omega $%
\textit{\ }$\backslash $\textit{\ }$[S_{\vec{F}}(u)$\textit{\ }$\cup $%
\textit{\ }$S_{\vec{F}}(v)].$

\endproof%

\bigskip

\proof
\textbf{(of Corollary B.1) }First we claim that both $S_{\vec{F}}(u)$ and $%
S_{\vec{F}}(v)$ are nowhere dense. Suppose the converse holds. Then we can
find a small ball $B$ contained in either $S_{\vec{F}}(u)$ or $S_{\vec{F}%
}(v),$ say $B$ $\subset $ $S_{\vec{F}}(u).$ That is, $\nabla u$ $+$ $\vec{F}$
$=$ $0$ in $B.$ It follows from $u_{IJ}$ $=$ $u_{JI}$ that $h_{IJ}$ :$=$ $%
\partial _{I}F_{J}-\partial _{J}F_{I}$ $=$ $0$ in $B$ for all $I,J.$\ So we
have $d\Theta _{u}$ $=$ $\frac{1}{2}h_{IJ}dx^{I}\wedge dx^{J}$ $=$ $0$ $=$ $%
d\Theta _{v}$ in $B,$\ contradicting the nonintegrability of $\Theta _{u}$
or $\Theta _{v}.$\ The above argument also works for $B$ $\subset $ $S_{\vec{%
F}}(v).$ So we have shown that both $S_{\vec{F}}(u)$ and $S_{\vec{F}}(v)$
are nowhere dense in $\Omega .$

Now we need only to show $\nabla u$ $=$ $\nabla v$ on $\Omega $\textit{\ }$%
\backslash $\textit{\ }$[S_{\vec{F}}(u)$\textit{\ }$\cup $\textit{\ }$S_{%
\vec{F}}(v)].$ We can invoke Theorem 5.1 in \cite{chy} to have $N_{\vec{F}%
}(u)$\textit{\ }$=$\textit{\ }$N_{\vec{F}}(v)$ in $\Omega $\textit{\ }$%
\backslash $\textit{\ }$[S_{\vec{F}}(u)$\textit{\ }$\cup $\textit{\ }$S_{%
\vec{F}}(v)]$ (comparing with the proof of Theorem A). By Theorem B we have $%
\nabla u$ $=$ $\nabla v$ in $\Omega $\textit{\ }$\backslash $\textit{\ }$[S_{%
\vec{F}}(u)$\textit{\ }$\cup $\textit{\ }$S_{\vec{F}}(v)].$

\endproof%

\bigskip

We remark that the condition $m$\textit{\ }$\geq $\textit{\ }$rank(h_{IJ})$%
\textit{\ }$\geq $\textit{\ }$3$ implies nonintegrability of $\Theta _{u}$
and $\Theta _{v}.$ Suppose one of them, say $\Theta _{u},$ is integrable.
Then (\ref{A10})\textit{\ }holds by Lemma 2.3$.$ It follows from Lemma 2.1
that $rank(h_{IJ}(p))$ $=$ $0$ or $2,$ a contradiction. Thus we have given
another proof of Theorem A by making use of Theorem B.

\bigskip

\proof
\textbf{(of Theorem A}$^{\prime },$ \textbf{Corollary A}$^{\prime }.1,$%
\textbf{\ Theorem B}$^{\prime },$\textbf{\ Corollary B}$^{\prime }.1$\textbf{%
)} Note that $|\nabla u+\vec{F}|$ $\geq $ $C$ and $|\nabla v+\vec{F}|$ $\geq 
$ $C$ for some constant $C$ $>$ $0$ imply%
\begin{equation*}
N_{\vec{F}}(u)=\frac{\nabla u+\vec{F}}{|\nabla u+\vec{F}|}\text{\textit{\
and }}N_{\vec{F}}(v)=\frac{\nabla v+\vec{F}}{|\nabla v+\vec{F}|}
\end{equation*}%
\noindent exist. Observe that $N_{\vec{F}}(u)$ $=$ $N_{\vec{F}}(v)$ in $%
\Omega $ (a.e.) implies $\nabla N_{\vec{F}}(u)$ $=$ $\nabla N_{\vec{F}}(v)$
in $\Omega $ (a.e.). Moreover, noting that $\nabla |\vec{X}|$ $=$ $\frac{%
\vec{X}}{|\vec{X}|}$ $\in $ $L^{\infty }$ where $\vec{X}$ $=$ $(x^{1},$ $%
x^{2},$ $...,$ $x^{m}),$ we have%
\begin{equation*}
\partial _{i}|\nabla u+\vec{F}|=\sum_{j=1}^{m}\frac{(\partial
_{j}u+F_{j})\partial _{i}(\partial _{j}u+F_{j})}{|\nabla u+\vec{F}|}\text{
(a.e.)}
\end{equation*}%
\noindent by Theorem 7.8 in \cite{GT}. Now apply the same reasoning as in
the proof of Theorem A (Corollary A.1, Theorem B, Corollary B.1, resp.) to
reach the conclusion.

\endproof%

\bigskip

\section{Proof of Theorem C}

The proof of Theorem C is similar as that of Theorem B in \cite{chy} as long
as we replace "$\ast "$ by "$b$" or $\vec{G}^{\ast }$ by $\vec{G}^{b}.$ We
recall the definition of $\vec{G}^{b}$ for $\vec{G}$ $=$ $(G_{1},$ $...,$ $%
G_{m})$ as follows:%
\begin{equation}
\vec{G}^{b}=(\sum_{k=1}^{m}a^{1k}G_{k},\sum_{k=1}^{m}a^{2k}G_{k},...,%
\sum_{k=1}^{m}a^{mk}G_{k})  \label{3.1}
\end{equation}

\noindent where $a^{jk\prime }s$ are real constants such that $a^{jk}+a^{kj} 
$ $=$ $0$ for $1$ $\leq $ $j,k$ $\leq $ $m.$ For the reader's convenience,
we will sketch the idea of the proof based on some reasonings in \cite{chy}.

Let $\kappa (\varepsilon )$ denote the Lebesgue measure of the set $S_{\vec{F%
}}(u_{\varepsilon })\cap \{\nabla \varphi \neq 0\}$ where $u_{\varepsilon }$ 
$=$ $u+\varepsilon \varphi ,$ $\varphi $ $=$ $v-u.$ There are at most
countably many $\varepsilon ^{\prime }s$ with $\kappa (\varepsilon )$ $>$ $0$
(see Section 3 in \cite{chy})$.$ We call such an $\varepsilon $ singular,
otherwise regular (i.e., $\kappa (\varepsilon )$ $=$ $0).$ Now we have

\bigskip

\textbf{Lemma 3.1 }(Theorem 5.1 in \cite{chy}) \textbf{\ }\textit{Let }$u,v$%
\textit{\ }$\in $\textit{\ }$W^{1,1}(\Omega )$\textit{\ be two minimizers
for }$\mathcal{F}_{H}(u)$\textit{\ such that }$u-v$\textit{\ }$\in $\textit{%
\ }$W_{0}^{1,1}(\Omega ).$\textit{\ Let }$u_{\varepsilon }$\textit{\ }$%
\equiv $\textit{\ }$u+\varepsilon (v-u).$\textit{\ Then for any pair of
regular }$\varepsilon _{1},$\textit{\ }$\varepsilon _{2}$\textit{\ }$\in $%
\textit{\ }$[0,1]$ \textit{(with respect to }$\varphi $\textit{\ }$=$\textit{%
\ }$v-u),$\textit{\ there holds }$N_{\vec{F}}(u_{\varepsilon _{1}})=N_{\vec{F%
}}(u_{\varepsilon _{2}})$\textit{\ in }$\Omega \backslash \lbrack S_{\vec{F}%
}(u_{\varepsilon _{1}})\cup S_{\vec{F}}(u_{\varepsilon _{2}})]$ $(a.e.).$

\bigskip

\textbf{Lemma 3.2. }\textit{Let }$u,v$\textit{\ }$\in $\textit{\ }$%
W^{1}(\Omega )$\textit{\ where the domain }$\Omega $\textit{\ is contained
in }$R^{m}.$\textit{\ Let }$u_{\varepsilon }$\textit{\ }$\equiv $\textit{\ }$%
u+\varepsilon (v-u).$\textit{\ Suppose }$N_{\vec{F}}(u_{\varepsilon
_{1}})=N_{\vec{F}}(u_{\varepsilon _{2}})$\textit{\ in }$\Omega \backslash
\lbrack S_{\vec{F}}(u_{\varepsilon _{1}})\cup S_{\vec{F}}(u_{\varepsilon
_{2}})]$\textit{\ for a pair }$\varepsilon _{1},$\textit{\ }$\varepsilon
_{2} $\textit{\ such that }$\varepsilon _{1}$\textit{\ }$\neq $\textit{\ }$%
\varepsilon _{2}$\textit{. Then for }$j=1,2,$\textit{\ there holds}

\begin{equation}
(\nabla u_{\varepsilon _{j}}+\vec{F})^{b}\cdot (\nabla v-\nabla u)=0\text{
in }\Omega \text{ (a.e.)}  \label{3.2}
\end{equation}

\bigskip

Lemma 3.2 extends Lemma 5.2 in \cite{chy}. Note that in deducing 
\begin{equation}
N_{\vec{F}}(u_{\varepsilon _{j}})^{b}\cdot \nabla u_{\varepsilon _{j}}=\frac{%
\vec{F}^{b}\cdot \nabla u_{\varepsilon _{j}}}{|\nabla u_{\varepsilon _{j}}+%
\vec{F}|}=\vec{F}^{b}\cdot N_{\vec{F}}(u_{\varepsilon _{j}})  \label{3.3}
\end{equation}

\noindent (cf. (5.4) in the proof of Lemma 5.2 in \cite{chy}), we have used
the property 
\begin{equation*}
\vec{G}^{b}\cdot \vec{G}=0
\end{equation*}%
\noindent twice$.$ $\vec{G}^{b}\cdot \vec{G}$ $=$ $0$ holds because%
\begin{eqnarray*}
\vec{G}^{b}\cdot \vec{G} &=&\sum_{j=1}^{m}a^{jk}G_{k}G_{j} \\
&=&\sum_{j=1}^{m}(-a^{kj})G_{k}G_{j}=-\vec{G}^{b}\cdot \vec{G}.
\end{eqnarray*}

Since $N_{\vec{F}}(u_{\varepsilon _{1}})=N_{\vec{F}}(u_{\varepsilon _{2}})$
in $\Omega \backslash \lbrack S_{\vec{F}}(u_{\varepsilon _{1}})\cup S_{\vec{F%
}}(u_{\varepsilon _{2}})]$ by assumption (hence $N_{\vec{F}}(u_{\varepsilon
_{1}})^{b}=N_{\vec{F}}(u_{\varepsilon _{2}})^{b}$ also), we take the
difference of (\ref{3.3}) for $j=1$ and $j=2$ to obtain%
\begin{equation}
N_{\vec{F}}(u_{\varepsilon _{1}})^{b}\cdot (\nabla u_{\varepsilon
_{2}}-\nabla u_{\varepsilon _{1}})=0.  \label{3.4}
\end{equation}

\noindent Formula (\ref{3.2}) for $j=1$ then follows from (\ref{3.4}) by
noting that $v-u$ $=$ $(u_{\varepsilon _{2}}-u_{\varepsilon
_{1}})/(\varepsilon _{2}-\varepsilon _{1}).$

\bigskip

\textbf{Lemma 3.3.} \textit{Let }$\Omega $\textit{\ be a bounded domain in }$%
R^{m}.$\textit{\ Let $w\in W_{0}^{1,p}(\Omega )$, $\sigma \in W^{1,q}(\Omega
)$, where $1\leq p<\infty $, $q=\frac{p}{p-1}$ ($q=\infty $ for $p=1$).} 
\textit{\ Let $\vec{F}$ (a vector field) $\in $ $W^{1,1}(\Omega )\cap
L^{q}(\Omega )$ satisfying $div\vec{F}^{b}$}\textit{\ }$>$\textit{\ }$0$%
\textit{\ (a.e.) or }$div\vec{F}^{b}$\textit{\ }$<$\textit{\ }$0$\ \textit{%
(a.e.).}\textit{\ Suppose }$(\nabla \sigma +\vec{F})^{b}\cdot \nabla w$%
\textit{\ }$=$\textit{\ }$0$\textit{\ in }$\Omega $\textit{\ (a.e.). Then }$%
w $\textit{\ }$\equiv $\textit{\ }$0$\textit{\ in }$\Omega $\textit{\ (a.e.).%
}

\bigskip

Replacing "*" by "$b$" in the proof of Theorem 5.3 in \cite{chy} gives a
proof of Lemma 3.3. We give an outline of the proof below. Approximate $w,$ $%
\sigma ,$ $\vec{F}$ by $\omega _{j}$ $\in $ $C_{0}^{\infty },$ $v_{k}$ $\in $
$C^{\infty },$ $\vec{F}_{\bar{k}}$ $\in $ $C^{\infty }$ in $W^{1,p},$ $%
W^{1,q},$ $W^{1,1}\cap L^{q},$ resp.. Suppose $\omega _{j}$ does not vanish
identically. Then for a decreasing sequence of $a_{i}$ $>$ $0$ converging to 
$0,$ $\Omega _{j,i}$ $:=$ $\{|\omega _{j}|$ $>$ $a_{i}\}$ $\subset \subset $ 
$\Omega $ is not empty for large $i$ and $\partial \Omega _{j,i}$ is $%
C^{\infty }$ smooth. Consider%
\begin{equation*}
I_{j,i,k,\bar{k}}:=\int_{\partial \Omega _{j,i}}|\omega _{j}|(\nabla v_{k}+%
\vec{F}_{\bar{k}})^{b}\cdot \nu
\end{equation*}%
\noindent where $\nu $ denotes the outer normal of $\partial \Omega _{j,i}.$
By using 
\begin{eqnarray*}
\func{div}(\nabla v_{k}+\vec{F}_{\bar{k}})^{b} &=&\func{div}(\nabla
v_{k})^{b}+\func{div}(\vec{F}_{\bar{k}})^{b} \\
&=&0+\func{div}(\vec{F}_{\bar{k}})^{b},
\end{eqnarray*}

\noindent we get%
\begin{equation*}
I_{j,i,k,\bar{k}}=a_{i}\int_{\Omega _{j,i}}\func{div}(\vec{F}_{\bar{k}})^{b}
\end{equation*}%
\noindent and hence%
\begin{equation}
\lim_{i\rightarrow \infty }I_{j,i,k,\bar{k}}=0  \label{3.5}
\end{equation}%
\noindent On the other hand, we compute 
\begin{eqnarray}
&&I_{j,i,k,\bar{k}}-\int_{\Omega \backslash \{\omega _{j}=0\}}\{\nabla
|\omega _{j}|\cdot (\nabla v_{k}+\vec{F}_{\bar{k}})^{b}+|\omega _{j}|\func{%
div}\vec{F}_{\bar{k}}^{b}\}  \label{3.6} \\
&=&(\int_{\Omega _{j,i}}-\int_{\Omega \backslash \{\omega _{j}=0\}})\{\nabla
|\omega _{j}|\cdot (\nabla v_{k}+\vec{F}_{\bar{k}})^{b}+|\omega _{j}|\func{%
div}\vec{F}_{\bar{k}}^{b}\}  \notag \\
&=&-\sum_{l=i}^{\infty }\int_{\Omega _{j,l+1}\backslash \Omega
_{j,l}}\{\nabla |\omega _{j}|\cdot (\nabla v_{k}+\vec{F}_{\bar{k}%
})^{b}+|\omega _{j}|\func{div}\vec{F}_{\bar{k}}^{b}\}  \notag \\
&=&-\sum_{l=i}^{\infty }(I_{j,l+1,k,\bar{k}}-I_{j,l.k,\bar{k}%
})=-\lim_{m\rightarrow \infty }I_{j,m,k,\bar{k}}+I_{j,i,k,\bar{k}}.  \notag
\end{eqnarray}

\noindent By (\ref{3.5}) and (\ref{3.6}) we have%
\begin{eqnarray}
0 &=&\int_{\Omega \backslash \{\omega _{j}=0\}}\{\nabla |\omega _{j}|\cdot
(\nabla v_{k}+\vec{F}_{\bar{k}})^{b}+|\omega _{j}|\func{div}\vec{F}_{\bar{k}%
}^{b}\}  \label{3.7} \\
&=&\int_{\Omega }\{\nabla |\omega _{j}|\cdot (\nabla v_{k}+\vec{F}_{\bar{k}%
})^{b}+|\omega _{j}|\func{div}\vec{F}_{\bar{k}}^{b}\}  \notag
\end{eqnarray}

\noindent in which we have used $\nabla |\omega _{j}|$ $=$ $0$ if $\omega
_{j}$ $=$ $0$ (p.152 in \cite{GT}). Letting $\bar{k}$ $\rightarrow $ $\infty
,$ $k$ $\rightarrow $ $\infty $ in (\ref{3.7}) gives%
\begin{equation*}
0=\int_{\Omega }\{\nabla |\omega _{j}|\cdot (\nabla \sigma +\vec{F}%
)^{b}+|\omega _{j}|\func{div}\vec{F}^{b}\}.
\end{equation*}

\noindent Since $(\nabla \sigma +\vec{F})^{b}\cdot \nabla w$\textit{\ }$=$%
\textit{\ }$0$ by assumption, we estimate%
\begin{eqnarray}
&&\int_{\Omega }\{\nabla |\omega _{j}|\cdot (\nabla \sigma +\vec{F})^{b}
\label{3.8} \\
&=&\int_{\{\omega _{j}>0\}}(\nabla \omega _{j}-\nabla w)\cdot (\nabla \sigma
+\vec{F})^{b}-\int_{\{\omega _{j}<0\}}(\nabla \omega _{j}-\nabla w)\cdot
(\nabla \sigma +\vec{F})^{b}  \notag \\
&\rightarrow &0\text{ \ as }j\rightarrow \infty .  \notag
\end{eqnarray}

\noindent On the other hand, we have%
\begin{equation}
\lim_{j\rightarrow \infty }\int_{\Omega }|\omega _{j}|\func{div}\vec{F}%
^{b}=\int_{\Omega }|w|\func{div}\vec{F}^{b}>0\text{ or }<0  \label{3.9}
\end{equation}

\noindent due to $\func{div}\vec{F}^{b}$ $>$ $0$ or $<$ $0$ (a.e.) by
assumption if $w$ $\neq $ $0.$ In view of (\ref{3.7}), (\ref{3.8}), and (\ref%
{3.9}), we reach a contradiction. Therefore $w$\textit{\ }$=$\textit{\ }$0$%
\textit{\ }in $\Omega $\ (a.e.).

\bigskip

\proof
\textbf{(of Theorem C)} The proof follows from Lemmas 3.1, 3.2, and 3.3 with 
$p$ $=$ $q$ $=$ $2,$ $\sigma $ $=$ $u_{\varepsilon _{1}},$ and $w$ $=$ $v-u.$

\endproof%

\bigskip

We would like to mention a result about the size of the singular set for $u$ 
$\in $ $C^{1},$ $\vec{F}$ $\in $ $C^{1}$ under the same condition on $\vec{F}
$ as in Theorem C.

\bigskip

\textbf{Proposition 3.4.} \textit{Let }$\Omega $\textit{\ be a bounded
domain of }$R^{m}.$\textit{\ Let }$u$\textit{\ }$\in $\textit{\ }$%
C^{1}(\Omega ),$\textit{\ }$\vec{F}$\textit{\ }$=$\textit{\ (}$F_{1},$%
\textit{\ }$F_{2},$\textit{\ }$...,$\textit{\ }$F_{m})$\textit{\ }$\in $%
\textit{\ }$C^{1}(\Omega ).$\textit{\ Suppose }%
\begin{equation*}
\func{div}\vec{F}^{b}=\sum_{j,k=1}^{m}a^{jk}\partial _{j}F_{k}>0\text{ ( }<0,%
\text{ resp.)}
\end{equation*}%
\textit{\noindent where }$a^{jk}$\textit{'s are real constants such that }$%
a^{jk}+a^{kj}=0.$\textit{\ Then }$S_{\vec{F}}(u)$\textit{\ is nowhere dense
in }$\Omega .$

\bigskip

\proof
Observe that $S_{\vec{F}}(u)$ is a closed set. So if $S_{\vec{F}}(u)$ is not
nowhere dense in $\Omega ,$ there there is a point $p_{1}$ $\in $ $S_{\vec{F}%
}(u)$ such that $S_{\vec{F}}(u)$ contains $B_{r_{1}}(p_{1}),$ a ball of
center $p_{1}$ with radius $r_{1}$ $>$ $0.$ Take a sequence of $C^{\infty }$
smooth functions $u_{k}$ converging to $u$ in $C^{1}$ norm on the closure of 
$B_{r_{2}}(p_{1})$ for $0$ $<$ $r_{2}$ $<$ $r_{1}.$ Let $\nu $ denote the
unit outer normal. Since $\nabla u$ $+$ $\vec{F}$ $=$ $0$ in $%
B_{r_{1}}(p_{1}),$ we have%
\begin{eqnarray*}
0 &=&\doint\limits_{\partial B_{r_{2}}(p_{1})}(\nabla u+\vec{F})^{b}\cdot \nu
\\
&=&\lim_{k\rightarrow \infty }\doint\limits_{\partial
B_{r_{2}}(p_{1})}(\nabla u_{k}+\vec{F})^{b}\cdot \nu \\
&=&\lim_{k\rightarrow \infty }\int_{B_{r_{2}}(p_{1})}\func{div}(\nabla u_{k}+%
\vec{F})^{b}\text{ (by the divergence theorem)} \\
&=&\int_{B_{r_{2}}(p_{1})}\func{div}\vec{F}^{b}\text{ (since }\func{div}%
(\nabla u_{k})^{b}=0) \\
&>&0\text{ (}<0,\text{ resp.),}
\end{eqnarray*}%
\textit{\noindent }a contradiction.

\endproof%

\bigskip

Note that Proposition 3.4 generalizes Lemma 3.1 in \cite{chmy2}.

\bigskip

\section{Proof of Theorem E and examples}

\proof
\textbf{(of Theorem E) }Let $U_{I}$ $=$ $D\nu _{I}-F_{I}.$ So $\nu _{I}$ $=$ 
$\frac{U_{I}+F_{I}}{D}.$ Since $\nu $ is a unit vector, we get $D$ $=$ ($%
\sum_{I=1}^{m}$($U_{I}+F_{I})^{2})^{1/2}.$ Let $\hat{U}_{I}$ $:=$ $%
U_{I}+F_{I}.$ By the same computation to reach (\ref{A5}), we have%
\begin{eqnarray}
&&\delta _{I}\nu _{J}-\delta _{J}\nu _{I}  \label{D1} \\
&=&\frac{\partial _{I}\hat{U}_{J}-\partial _{J}\hat{U}_{I}}{D}-\frac{\hat{U}%
_{J}\hat{U}_{K}(\partial _{I}\hat{U}_{K}-\partial _{K}\hat{U}_{I})}{D^{3}} 
\notag \\
&&-\frac{\hat{U}_{I}\hat{U}_{K}(\partial _{K}\hat{U}_{J}-\partial _{J}\hat{U}%
_{K})}{D^{3}}.  \notag
\end{eqnarray}

\noindent Noting that $\frac{\hat{U}_{J}}{D}$ $=$ $\nu _{J}$ and
substituting $\hat{U}_{I}$ $:=$ $U_{I}+F_{I}$ into (\ref{D1}), we obtain%
\begin{eqnarray}
&&\delta _{I}\nu _{J}-\delta _{J}\nu _{I}  \label{D2} \\
&=&\frac{U_{IJ}-\nu _{J}\nu _{K}U_{IK}-\nu _{I}\nu _{K}U_{KJ}}{D}  \notag \\
&&+\frac{h_{IJ}-\nu _{J}\nu _{K}h_{IK}-\nu _{I}\nu _{K}h_{KJ}}{D}  \notag
\end{eqnarray}

\noindent where $U_{IJ}$ $:=$ $\partial _{I}U_{J}-\partial _{J}U_{I}$ and
recall $h_{IJ}$ $=$ $\partial _{I}F_{J}-\partial _{J}F_{I}.$ By the
assumption (\ref{1.1.4}) and (\ref{D2}), we have%
\begin{equation}
U_{IJ}-\nu _{J}\nu _{K}U_{IK}-\nu _{I}\nu _{K}U_{KJ}=0  \label{D3}
\end{equation}

Let $U$ $=$ $(U_{IJ}).$ Recall that we view $\nu $ $=$ $(\nu _{J})$ as a ($%
m\times 1$) unit column real vector and $\nu ^{T},$ the transpose of $\nu ,$
as a ($1\times m)$ unit row vector. We can write (\ref{D3}) as follows:%
\begin{equation}
U=(U\nu )\nu ^{T}-\nu (U\nu )^{T}  \label{D3.1}
\end{equation}

\noindent in which we have used skew-symmetry of $U$ (i.e., $U^{T}$ $=$ $-U$
where $U^{T}$ denotes the transpose of $U).$ In terms of differential forms,
we have%
\begin{eqnarray}
d(U_{I}dx^{I}) &=&\frac{1}{2}U_{IJ}dx^{I}\wedge dx^{J}  \label{D3.2} \\
&=&\frac{1}{2}\{(U\nu )_{I}\nu _{J}-\nu _{I}(U\nu )_{J}\}dx^{I}\wedge dx^{J}
\notag \\
&=&\frac{1}{2}\{(U\nu )^{\#}\wedge \nu ^{\#}-\nu ^{\#}\wedge (U\nu )^{\#}\} 
\notag \\
&=&(U\nu )^{\#}\wedge \nu ^{\#}  \notag
\end{eqnarray}

\noindent by (\ref{D3.1}). Recall\textit{\ }that $w^{\#}$ denotes the
corresponding 1-form for a vector $w.$

Now we observe that%
\begin{equation*}
d(D\nu ^{\#})=d(U_{I}dx^{I})+dF^{\#}.
\end{equation*}%
\textit{\noindent }Comparing with condition (\ref{1.1.5}), we have%
\begin{equation}
\nu ^{\#}\lrcorner d(U_{I}dx^{I})=0.  \label{D5}
\end{equation}%
\textit{\noindent }Substituting (\ref{D3.2}) into (\ref{D5}), we have%
\begin{eqnarray}
0 &=&\nu ^{\#}\lrcorner ((U\nu )^{\#}\wedge \nu ^{\#})  \label{D5.1} \\
&=&<\nu ^{\#},(U\nu )^{\#}>\nu ^{\#}-<\nu ^{\#},\nu ^{\#}>(U\nu )^{\#} 
\notag \\
&=&-(U\nu )^{\#}.  \notag
\end{eqnarray}
\textit{\noindent }Here we have used $<\nu ^{\#},\nu ^{\#}>$ $=$ $<\nu ,\nu
> $ $=$ $1$ and $<\nu ^{\#},(U\nu )^{\#}>$ $=$ $<\nu ,U\nu >$ $=$ $0$ since%
\begin{eqnarray}
&<&\nu ,U\nu >=<U^{T}\nu ,\nu >  \label{D5.2} \\
&=&<-U\nu ,\nu >=-<\nu ,U\nu >.  \notag
\end{eqnarray}%
\textit{\noindent }It follows from (\ref{D5.1}) that $U\nu $ $=$ $0,$ and
hence $U$ $=$ $0$ by (\ref{D3.1}), i.e., $U_{IJ}$ $=$ $0.$ This means $%
\partial _{I}U_{J}$ $=$ $\partial _{J}U_{I}.$ Therefore locally we can find
a ($C^{\infty }$ smooth) function $u$ such that $U_{I}$ $=$ $\partial _{I}u.$
By the definition of $U_{I},$ we have $\nu $\textit{\ }$=$\textit{\ }$\frac{%
\nabla u+\vec{F}}{D}.$ We have completed the proof of Theorem E.

\endproof%

\bigskip

\textbf{Proposition 4.1}. \textit{Let }$\nu $\textit{\ be a unit vector.} 
\textit{Condition (\ref{1.1.5}) in Theorem E is equivalent to the following
system of first order equations in }$D:$%
\begin{equation}
\delta _{K}D=\nu _{J}(\partial _{J}\nu _{K})D-\nu _{J}h_{JK}  \label{D6}
\end{equation}%
\textit{\noindent for }$1$\textit{\ }$\leq $\textit{\ }$K$\textit{\ }$\leq $%
\textit{\ }$m.$

\bigskip

\proof
Observe that 
\begin{eqnarray}
&&\nu ^{\#}\lrcorner d(D\nu ^{\#}-F^{\#})  \label{D7} \\
&=&<\nu ^{\#},dD>\nu ^{\#}-dD<\nu ^{\#},\nu ^{\#}>+D\nu ^{\#}\lrcorner d\nu
^{\#}-\nu ^{\#}\lrcorner dF^{\#}.  \notag
\end{eqnarray}%
\textit{\noindent }Let $D_{I}$ denote $\partial _{I}D.$ We compute%
\begin{eqnarray}
&<&\nu ^{\#},dD>\nu ^{\#}  \label{D8} \\
&=&\nu _{I}D_{I}\nu _{K}dx^{K},  \notag
\end{eqnarray}

\begin{eqnarray}
\nu ^{\#}\lrcorner d\nu ^{\#} &=&\nu _{I}dx^{I}\lrcorner ((\partial _{K}\nu
_{J})dx^{K}\wedge dx^{J})  \label{D9} \\
&=&\nu _{I}(\partial _{K}\nu _{J})(\delta ^{IK}dx^{J}-\delta ^{IJ}dx^{K}) 
\notag \\
&=&\nu _{J}(\partial _{J}\nu _{K}-\partial _{K}\nu _{J})dx^{K},  \notag
\end{eqnarray}

\begin{eqnarray}
\nu ^{\#}\lrcorner dF^{\#} &=&\nu _{I}dx^{I}\lrcorner ((\partial
_{J}F_{K})dx^{J}\wedge dx^{K})  \label{D10} \\
&=&\nu _{I}(\partial _{J}F_{K})(\delta ^{IJ}dx^{K}-\delta ^{IK}dx^{J}) 
\notag \\
&=&\nu _{I}h_{IK}dx^{K}  \notag
\end{eqnarray}%
\textit{\noindent }in which we recall that $\delta ^{IJ}$ denotes the
Kronecker delta and $h_{IK}$ $:=$ $\partial _{I}F_{K}-\partial _{K}F_{I}.$
Substituting (\ref{D8}), (\ref{D9}), and (\ref{D10}) into (\ref{D7}), we
reduce (\ref{1.1.5}) to the following equations:%
\begin{eqnarray}
0 &=&\nu _{K}\nu _{I}D_{I}-D_{K}  \label{D11} \\
&&+\nu _{J}(\partial _{J}\nu _{K}-\partial _{K}\nu _{J})D-\nu _{I}h_{IK} 
\notag
\end{eqnarray}%
\textit{\noindent }for\textit{\ }$1$\textit{\ }$\leq $\textit{\ }$K$\textit{%
\ }$\leq $\textit{\ }$m.$ Noting that $\delta _{K}$ $:=$ $\partial _{K}-\nu
_{K}\nu _{I}\partial _{I}$ and $\nu _{J}\partial _{K}\nu _{J}$ $=$ $0$ due
to $\sum_{J}\nu _{J}^{2}$ $=$ $1$ in (\ref{D11}), we get (\ref{D6}).

\endproof%

\bigskip

Let us discuss the ($p$-area) situation of dimension 2 for $\vec{F}$ $=$ $%
(-y,$ $x)$. Write $\nu _{1}$ $=$ $\cos \theta ,$ $\nu _{2}$ $=$ $\sin \theta 
$ and $\nu _{\perp }$ $=$ $\nu _{2}\partial _{x}$ $-$ $\nu _{1}\partial _{y}$
$=$ $\sin \theta \partial _{x}$ $-$ $\cos \theta \partial _{y}.$ We have $%
d(D\nu ^{\#})$ $=$ $d(D\cos \theta dx+D\sin \theta dy)$ $=$ [$\partial
_{x}(D\sin \theta )$ $-$ $\partial _{y}(D\cos \theta )]$ $dx\wedge dy,$ $%
dF^{\#}$ $=$ $2dx\wedge dy,$ and hence%
\begin{eqnarray}
&&\nu ^{\#}\lrcorner d(D\nu ^{\#}-F^{\#})  \label{D11.1} \\
&=&\nu ^{\#}\lrcorner \{\func{div}(D\nu _{\perp })-2\}dx\wedge dy  \notag \\
&=&-\{\func{div}(D\nu _{\perp })-2\}\nu _{\perp }^{\#}.  \notag
\end{eqnarray}%
\textit{\noindent }It follows that $\nu ^{\#}\lrcorner d(D\nu ^{\#}-F^{\#})$ 
$=$ $0$ if and only if $\func{div}(D\nu _{\perp })-2$ $=$ $0$ if and only if 
$d(D\nu ^{\#}-F^{\#})$ $=$ $0.$ The equation $\func{div}(D\nu _{\perp })-2$ $%
=$ $0$ is in fact a basic (Codazzi-like) equation for a surface in $3$%
-dimensional Heisenberg group (see, e.g., (1.17) in \cite{chmy2} where $V$
is supposed to be $\nu _{\perp }$ here). In higher dimensions, we can have
examples satisfying (\ref{1.1.5}), but $d(D\nu ^{\#}-F^{\#})$ $\neq $ $0.$

\bigskip

\textbf{Example 4.2}. In dimension $m$ $=$ $4$ we take $\vec{F}$ $=$ $%
(-y^{1},$ $x^{1},$ $-y^{2},$ $x^{2})$ ($p$-area situation). Let $\nu ^{\#}$ $%
=$ $dx^{1}$ and $D$ $=$ $-2y^{1}$ ($>$ $0$ in the region of $y^{1}$ $<$ $0).$
We compute%
\begin{eqnarray*}
&&d(D\nu ^{\#}-F^{\#}) \\
&=&-2dy^{1}\wedge dx^{1}-2(dx^{1}\wedge dy^{1}+dx^{2}\wedge dy^{2}) \\
&=&-2dx^{2}\wedge dy^{2}\neq 0.
\end{eqnarray*}%
\textit{\noindent }Clearly $\nu ^{\#}\lrcorner d(D\nu ^{\#}-F^{\#})$ $=$ $%
dx^{1}\lrcorner (-2dx^{2}\wedge dy^{2})$ $=$ $0.$

\bigskip

\textbf{Example 4.3}. We can have examples satisfying (\ref{1.1.4}), but not
(\ref{1.1.5}). Take $\vec{F}$ $=$ $0$ and $\nu $ $=$ a constant unit vector$%
. $ So we have $\delta _{I}\nu _{J}-\delta _{J}\nu _{I}$ $=$ $0$ while $%
h_{IJ}$ $=$ $\partial _{I}F_{J}-\partial _{J}F_{I}$ $=$ $0.$ Therefore (\ref%
{1.1.4}) holds. Choose $\nu $ (constant unit) such that we can pick up
another unit vector $\nu _{\perp }$ perpendicular to $\nu $ with the
property: ($\nu _{\perp })_{J}$ $>$ $0$ for all $J.$ Take $D$ $=$ ($\nu
_{\perp })_{J}x^{J}$ (summation convention) $>$ $0$ in the region of all $%
x^{J}$ $>$ $0.$ It follows that $dD$ $=$ ($\nu _{\perp })_{J}dx^{J}$ $=$ $%
\nu _{\perp }^{\#},$ and hence 
\begin{eqnarray}
&<&\nu ^{\#},dD>  \label{D12} \\
&=&<\nu ^{\#},\nu _{\perp }^{\#}>=<\nu ,\nu _{\perp }>=0.  \notag
\end{eqnarray}

We\textit{\ }can now compute%
\begin{eqnarray*}
&&\nu ^{\#}\lrcorner d(D\nu ^{\#}-F^{\#}) \\
&=&\nu ^{\#}\lrcorner d(D\nu ^{\#})\text{ \ \ }(\because \vec{F}=0) \\
&=&<\nu ^{\#},dD>\nu ^{\#}-<\nu ^{\#},\nu ^{\#}>dD\text{ }(\because d\nu
^{\#}=0) \\
&=&0-dD=-\nu _{\perp }^{\#}\neq 0\text{ }(\text{by (\ref{D12})).}
\end{eqnarray*}%
\textit{\noindent }I.e.\textit{, }(\ref{1.1.5}) does not hold. Note that for
such ($\nu ,D,\vec{F}),$ $\nu $ $\neq $ $\frac{\nabla u+\vec{F}}{D}.$for any
function $u$ (recall that (\ref{1.1.5}) is a necessary condition.for $\nu $ $%
=$ $\frac{\nabla u+\vec{F}}{D}$ for some $u).$

\bigskip

\section{Appendix}

In this section we collect some more facts about the properties of $U$
satisfying (\ref{D3.1}) (or (\ref{D3})). Recall that the rank of a matrix $%
U, $ denoted as $rank(U),$ is the dimension of the range $Range(U)$ (or
image) of $U.$ Let $||w||$ $=$ $<w,w>^{1/2}.$

\bigskip

\textbf{Proposition A.1}. \textit{Let }$U$\textit{\ be an }$m\times m$%
\textit{\ real matrix (}$\mathit{m\geq 2)}$ \textit{such that }$U=-U^{T}$%
\textit{\ (skew-symmetric) and }$rank(U)$\textit{\ }$=$\textit{\ }$2.$ 
\textit{Then }$U\nu $\textit{\ }$\neq $\textit{\ }$0$\textit{\ for some }$%
\nu $\textit{\ }$\neq $\textit{\ }$0$\textit{\ and for such }$\nu ,$\textit{%
\ we have}

\textit{(1) }$U^{2}\nu \neq 0;$

\textit{(2) }$<\nu ,U\nu >$\textit{\ }$=$\textit{\ }$<U\nu ,U^{2}\nu >$%
\textit{\ }$=$\textit{\ }$0;$

\textit{(3) }$Range(U)$\textit{\ is spanned by }$U\nu $\textit{\ and }$%
U^{2}\nu ;$

\textit{(4) }$Range(U^{2})$\textit{\ is also spanned by }$U\nu $\textit{\
and }$U^{2}\nu ,$\textit{\ in particular, }$rank(U^{2})$\textit{\ }$=$%
\textit{\ }$2;$

\textit{(5) }$U\nu $\textit{\ and }$U^{2}\nu $\textit{\ are eigenvectors of }%
$U^{2}$\textit{\ with the same eigenvalue}

\begin{equation}
\mathit{\rho =-}\frac{||U^{2}\nu ||^{2}}{||U\nu ||^{2}}\mathit{.}
\label{A-1}
\end{equation}

\bigskip

\proof
If $U^{2}\nu $ $=$ $0,$ then%
\begin{eqnarray*}
0 &=&<U^{2}\nu ,\nu > \\
&=&<U\nu ,U^{T}\nu >=-<U\nu ,U\nu >.
\end{eqnarray*}%
\textit{\noindent }So $U\nu $ $=$ $0,$ a contradiction. We have proved (1).
Since $U$ is skew-symmetric, we have%
\begin{eqnarray*}
&<&w,Uw>=<U^{T}w,w> \\
&=&-<Uw,w>=-<w,Uw>,\text{ and hence}
\end{eqnarray*}%
\begin{equation}
<w,Uw>=0  \label{A0}
\end{equation}%
\textit{\noindent }for\textit{\ }any $w.$ Substituting $w$ $=$ $\nu $ and $%
U\nu ,$ resp. in (\ref{A0})$,$ we get (2). By (1) and (2), $U\nu $ and $%
U^{2}\nu $ form an orthogonal basis for $Range(U).$ (3) follows. Next $%
U^{3}\nu $ $\neq $ $0$ by a similar argument in deducing (1). By (\ref{A0})
with $w$ $=$ $U^{2}\nu ,$ we get%
\begin{equation}
<U^{2}\nu ,U^{3}\nu >=0.  \label{A0.1}
\end{equation}%
\textit{\noindent }It follows that $U^{2}\nu ,$ $U^{3}\nu $ (=$U^{2}(U\nu ))$
are independent nonzero elements in $Range(U^{2}).$ On the other hand,
observe that $Range(U^{2})$ $\subset $ $Range(U),$ and hence $rank(U^{2})$ $%
\leq $ $2.$ Therefore $rank(U^{2})$ $=$ $2$ and $Range(U^{2})$ $=$ $Range(U)$
is also spanned by $U\nu $\textit{\ }and $U^{2}\nu $ by (3). We have proved
(4). Since $U^{3}\nu $ $\in $ $Range(U^{2})$ is perpendicular to $U^{2}\nu $
by (\ref{A0.1})$,$ we conclude that 
\begin{equation}
U^{2}(U\nu )=U^{3}\nu =\rho U\nu  \label{A0.11}
\end{equation}%
\textit{\noindent }for some $\rho $ $\in $ $R.$ It follows that%
\begin{eqnarray}
\rho &<&U\nu ,U\nu >=<U^{3}\nu ,U\nu >  \label{A0.2} \\
&=&<U^{2}\nu ,U^{T}U\nu >  \notag \\
&=&-<U^{2}\nu ,U^{2}\nu >.  \notag
\end{eqnarray}%
\textit{\noindent }Observe that $U^{2}(U^{2}\nu )$ $=$ $U(U^{3}\nu )$ $=$ $%
U(\rho U\nu )$ $=$ $\rho U^{2}\nu $ by (\ref{A0.11}). So $U\nu $\ and $%
U^{2}\nu $\ are eigenvectors of $U^{2}$ with the same eigenvalue $\rho .$
Formula (\ref{A-1}) follows from (\ref{A0.2}). We have proved (5).

\endproof%

\bigskip

\textbf{Proposition A.2.} \textit{Let }$U$\textit{\ be a nonzero
skew-symmetric real }$m\times m$\textit{\ matrix (}$m\geq 2),$\textit{\
i.e., }$U=-U^{T}$ \textit{and} $U\neq 0.$\textit{\ Then }$rank(U)$\textit{\ }%
$=$\textit{\ }$2$ \textit{if and only if}%
\begin{equation}
U=(U\nu )\nu ^{T}-\nu (U\nu )^{T}  \label{A1}
\end{equation}%
\textit{\noindent for any (}$m\times 1$\textit{) unit column real vector }$%
\nu $ \textit{satisfying}%
\begin{equation}
U^{2}\nu =\rho \nu  \label{A1-1}
\end{equation}%
\textit{\noindent for a nonzero real number }$\rho .$

\bigskip

\proof
Suppose $rank(U)$\textit{\ }$=$\textit{\ }$2.$ Then $Uw$ $\neq $ $0$ for
some $w$ $\neq $ $0.$ Take 
\begin{equation*}
\nu =\frac{Uw}{||Uw||}.
\end{equation*}%
\textit{\noindent }By Proposition A.1 (5) (with $\nu $ replaced by $w$ there$%
),$ we learn that $\nu $ and $U\nu $ are eigenvectors of $U^{2}$ with
nonzero eigenvalue $\rho $ (so (\ref{A1-1}) holds) and moreover,

\begin{equation}
0\neq \rho =-\frac{||U\nu ||^{2}}{||\nu ||^{2}}=-||U\nu ||^{2}.  \label{A1.0}
\end{equation}%
By Proposition A.1 (4) and (5), we learn that $0$ is the only eigenvalue
different from $\rho $ and the dimension of its eigenspace is $m-2.$ Let $%
\nu _{j},$ $j$ $=$ $3,$ $...,$ $m,$ be orthonormal eigenvectors of $U^{2}$
with eigenvalue $0.$ By Proposition A.1 (1)$,$ we have $U\nu _{j}$ $=$ $0$
(otherwise, $U^{2}\nu _{j}$ $\neq $ $0).$ Let 
\begin{equation*}
\tilde{U}=(U\nu )\nu ^{T}-\nu (U\nu )^{T}.
\end{equation*}%
\textit{\noindent }It\textit{\ }is now a direct verification that $\tilde{U}%
\nu _{j}$ $=$ $0,$ $j$ $=$ $3,$ $...,$ $m,$ since $<\nu ,\nu _{j}>$ $=$ $0$
and $<U\nu ,\nu _{j}>$ $=$ $<\nu ,U^{T}\nu _{j}>$ $=$ $-<\nu ,U\nu _{j}>$ $=$
$0$ by $U\nu _{j}$ $=$ $0.$ On the other hand, we have%
\begin{eqnarray*}
\tilde{U}\nu &=&(U\nu )<\nu ,\nu >-\nu <U\nu ,\nu > \\
&=&U\nu
\end{eqnarray*}%
\textit{\noindent }by $<\nu ,\nu >$ $=$ $1$ and $<U\nu ,\nu >$ $=$ $0.$ We
also compute%
\begin{eqnarray*}
&&\tilde{U}(U\nu ) \\
&=&U\nu <\nu ,U\nu >-\nu <U\nu ,U\nu > \\
&=&0+U(U\nu ).
\end{eqnarray*}%
\textit{\noindent }In the last equality, we have used (\ref{A0}), (\ref{A1.0}%
), and (\ref{A1-1})$.$ Altogether we conclude that $\tilde{U}$ $=$ $U.$ We
have shown (\ref{A1}). The reverse direction is due to Lemma 2.1.

\endproof%

\bigskip

Note that Proposition A.2 includes the converse of Lemma 2.1. In the
following Proposition we point out that (\ref{A1-1}) with $\rho $ given by (%
\ref{A1.0}) is also a necessary condition for (\ref{A1}) to hold. Note that
equation (\ref{A1}) is equivalent to%
\begin{equation}
U-U\nu \nu ^{T}-\nu \nu ^{T}U=0.  \label{A1.1}
\end{equation}

Let $\nu _{\perp }$ $=$ $\frac{U\nu }{||U\nu ||}.$ It follows from
skew-symmetry of $U$ that $<\nu _{\perp },\nu >$ $=$ $0.$

\bigskip

\textbf{Proposition A.3.} \textit{Let }$U$\textit{\ be a skew-symmetric real 
}$m\times m$\textit{\ matrix (}$m\geq 2)$\textit{\ such that (\ref{A1}) (or (%
\ref{A1.1})) holds. Then }

\textit{(1) }$U^{2}\nu =-||U\nu ||^{2}\nu $\textit{;}

\textit{(2) }$U$\textit{\ }$=$\textit{\ }$||U\nu ||(\nu _{\perp }\nu
^{T}-\nu (\nu _{\perp })^{T})$\textit{.}

\bigskip

\proof
Apply (\ref{A1}) to $U\nu $ to get $U^{2}\nu $ $=$ $U\nu <\nu ,U\nu >$ $-$ $%
\nu <U\nu ,U\nu >$ $=$ $0$ $-||U\nu ||^{2}\nu $. (1) follows. Substituting $%
U\nu $ $=$ $||U\nu ||\nu _{\perp }$ into (\ref{A1}) gives (2).

\endproof%


\bigskip

\begin{thebibliography}{99}
\bibitem{ASCV} L. Ambrosio, F. Serra Cassano, and D. Vittone, Intrinsic
regular hypersurfaces in Heisenberg groups, J. Geom. Anal. 16 (2006) 187-232.

\bibitem{Ba} Z. M. Balogh, Size of characteristic sets and functions with
prescribed gradient, J. reine angew. Math., 564 (2003) 63-83.

\bibitem{BSC} F. Bigolin and F. Serra Cassano, Distributional solutions of
Burgers' equation and intrinsic regular graphs in Heisenberg groups, J.
Math. Anal. Appl. 366 (2010) 561-568.

\bibitem{CDG} L. Capogna, D. Danielli, N. Garofalo, The geometric Sobolev
embedding for vector fields and the isoperimetric inequality, Comm. Anal.
Geom., 2 (1994) 203-215.

\bibitem{ch} J.-H. Cheng and J.-F. Hwang, Properly embedded and immersed
minimal surfaces in the Heisenberg group, Bull. Aus. Math. Soc., 70 (2004)
507-520.

\bibitem{ch1} J.-H. Cheng and J.-F. Hwang, Variations of generalized area
functionals and p-area minimizers of bounded variation in the Heisenberg
group, Bulletin of the Institute of Mathematics, Academia Sinica, New
Series, 5 (2010) 369-412.

\bibitem{chmy} J.-H. Cheng, J.-F. Hwang, A. Malchiodi, and P. Yang, Minimal
surfaces in pseudohermitian geometry, Ann Scuola Norm. Sup. Pisa Cl. Sci.
(5) Vol. IV (2005) 129-177.

\bibitem{chmy2} J.-H. Cheng, J.-F. Hwang, A. Malchiodi, and P. Yang, A
Codazzi-like equation and the singular set for $C^{1}$ smooth surfaces in
the Heisenberg group, to appear in J. reine angew. Math., 2012.

\bibitem{chy} J.-H. Cheng, J.-F. Hwang, and P. Yang, Existence and
uniqueness for p-area minimizers in the Heisenberg group, Math. Annalen, 337
(2007) 253-293.

\bibitem{chy1} J.-H. Cheng, J.-F. Hwang, and P. Yang, Regularity of $C^{1}$
smooth surfaces with prescribed $p$-mean curvature in the Heisenberg group,
Math. Annalen, 344 (2009) 1-35.

\bibitem{Chiu} H.-L. Chiu and S.-H. Lai, The fundamental theorems for curves
and surfaces in 3D Heisenberg group, arXiv:1301.6463.

\bibitem{DGN} D. Danielli, N. Garofalo, and D.-M. Nhieu, Minimal surfaces,
surfaces of constant mean curvature and isoperimetry in Carnot groups,
preprint, 2001.

\bibitem{FSS} B. Franchi, R. Serapioni, and F. Serra Cassano, Rectifiability
and perimeter in the Heisenberg group, Math. Ann., 321 (2001) 479-531.

\bibitem{GN} N. Garofalo and D.-M. Nhieu, Isoperimetric and Sobolev
inequalities for Carnot- Caratheodory spaces and the existence of minimal
surfaces, Comm. Pure Appl. Math., 49 (1996) 1081-1144.

\bibitem{GT} D. Gilbarg and N. S. Trudinger, Elliptic partial differential
equations of second order, 2nd ed., G.M.W. 224, Springer-Verlag, 1983.

\bibitem{LM} G. Leonardi and S. Masnou, On the isoperimetric problem in the
Heisenberg group $H^{n}$, Ann. Mat. Pura Appl. (4), No.4, 184 (2005) 533-553.

\bibitem{LR} G. Leonardi and S. Rigot, Isoperimetric sets on Carnot groups,
Houston J. Math., No.3, 29 (2003) 609-637.

\bibitem{MM} U. Massari and M. Miranda, Minimal surfaces of codimension one,
North-Holland, Amsterdam-New York-Oxford, 1984.

\bibitem{MR} R. Monti and M. Rickly, Convex isoperimetric sets in the
Heisenberg group, Ann. Sc. Norm. Super. Pisa Cl. Sci. (5), No.2, 8 (2009)
391-415.

\bibitem{Pan} P. Pansu,\emph{\ }Une inegalite isoperimetrique sur le groupe
de Heisenberg, C. R. Acad. Sci. Paris S\'{e}r. I Math., No.2, 295 (1982)
127-130.

\bibitem{Pau} S. D. Pauls, Minimal surfaces in the Heisenberg group,
Geometric Dedicata, 104 (2004) 201-231.

\bibitem{PCTV} A. Pinamonti, F. Serra Cassano, G. Treu, and D. Vittone, BV
minimizers of the area functional in the Heisenberg group under the bounded
slope condition, preprint.

\bibitem{RR} M. Ritor\'{e} and C. Rosales, Area-stationary surfaces in the
Heisenberg group $H^{1},$ Advances in Math., 219 (2008) 633-671.

\bibitem{SCV} F. Serra Cassano and D. Vittone, Graphs of bounded variation,
existence and local boundedness of nonparametric minimal surfaces in
Heisenberg groups, preprint.
\end{thebibliography}
\end{document}